\documentclass{conm-p-l}
\usepackage{amssymb}

\usepackage{amscd}
\usepackage{xy}

\xyoption{arrow}
\xyoption{matrix}
\xyoption{ps}
\xyoption{dvips}

\nonstopmode

\numberwithin{equation}{section}
\newtheorem{theorem}{Theorem}[section]
\newtheorem{lemma}[theorem]{Lemma}
\newtheorem{proposition}[theorem]{Proposition}
\newtheorem{corollary}[theorem]{Corollary}

\theoremstyle{definition}
\newtheorem{definition}[theorem]{Definition}
\newtheorem{remark}[theorem]{Remark}
\newtheorem{example}[theorem]{Example}

\newcommand\codim{\operatorname{codim}}
\newcommand\coker{\operatorname{coker}}
\newcommand\Ann{\operatorname{Ann}}

\newcommand\Hom{\operatorname{Hom}}
\newcommand\Ext{\operatorname{Ext}}
\newcommand\reg{\operatorname{reg}}

\newcommand\rank{\operatorname{rank}}
\newcommand\rankk{\operatorname{rank}_K}
\newcommand\depth{\operatorname{depth}}
\newcommand\soc{\operatorname{soc}}
\newcommand{\xx}{\underline x}

\newcommand\im{\operatorname{im}}

\newcommand{\HH}{H_{\mathfrak m}}
\newcommand{\hh}{h_{\mathfrak m}}

\newcommand{\TR}{\operatorname{Tor}^R}
\newcommand{\ER}{\operatorname{Ext}_R}
\newcommand{\Proj}{\operatorname{Proj}}

\newcommand{\s}{\; | \;}
\newcommand{\mif}{\mbox{if} ~}

\newcommand{\vep}{\varepsilon}
\newcommand{\ffi}{\varphi}

\newcommand{\cJ}{{\mathcal J}}

\newcommand{\cE}{{\mathcal E}}
\newcommand{\cF}{{\mathcal F}}

\newcommand{\cL}{{\mathcal L}}

\newcommand{\cO}{{\mathcal O}}
\newcommand{\cP}{{\mathcal P}}

\newcommand{\fm}{{\mathfrak m}}
\newcommand{\fa}{{\mathfrak a}}
\newcommand{\fc}{{\mathfrak c}}

\newcommand{\cOP}{\Omega_{\mathbb{P}^n}}

\newcommand {\ZZ}{\mathbb{Z}}

\newcommand {\PP}{\mathbb{P}}

\begin{document}
\title[Characterization by locally free resolutions]{Characterization
of some
  projective subschemes by locally free resolutions}
\author[Uwe Nagel]{Uwe Nagel}

\address{Fachbereich Mathematik und Informatik,
Universit\"at Paderborn,  D--33095 Paderborn, Germany}
\email{uwen@uni-paderborn.de}

\address{(Address since August 2002:
Department of Mathematics, University of Kentucky, Lexington, KY 40506-0027,
USA) }





\begin{abstract} A locally free resolution of a subscheme is by definition
  an exact sequence consisting of locally free sheaves (except the ideal
  sheaf) which has uniqueness properties like a free resolution. The
  purpose of this paper is to characterize certain locally Cohen-Macaulay
  subschemes by means of locally free resolutions. First we achieve this for
arithmetically Buchsbaum  subschemes. This leads to the notion of an
$\Omega$-resolution and extends the main  result of Chang in
\cite{Chang-Diff-Geom}. Second we characterize quasi-Buchsbaum subschemes
by means of weak  $\Omega$-resolutions. Finally, we describe the weak
$\Omega$-resolutions which belong to arithmetically Buchsbaum surfaces of
codimension two. Various applications of our results are given.
\end{abstract}


\maketitle

\dedicatory{Dedicated to the memory of Wolfgang Vogel}


\tableofcontents

\section{Introduction} \label{intro}

Since  free
resolutions have been invented by Hilbert they have proved to be very
useful in 
algebraic geometry and commutative algebra. Indeed, the  minimal free
resolution of a subscheme reflects its geometric properties
 (cf., for example, \cite{G}, \cite{Green-L-normality}, \cite{ADHPR-1997})
 as well as its algebraic properties. For example, the
 minimal free resolution can be used to characterize certain arithmetically
 Cohen-Macaulay subschemes of projective space. In this paper we want to
 establish some characterization of non-arithmetically
 Cohen-Macaulay subschemes by considering more general resolutions.

Let $\PP^n$ denote the projective space over a field $K$ and
let $X \subset \PP^n$ be a projective subscheme of codimension $c$ with a
minimal free resolution
$$
0 \to \cF_s \to \cF_{s-1} \to \ldots \to \cF_1 \to \cJ_X \to 0
$$
where $\cJ_X$ is the ideal sheaf of $X \subset \PP^n$ and the sheaves
$\cF_i$ are direct sums of line bundles on $\PP^n$. Due to the
Auslander-Buchsbaum formula we know
that it holds $s \geq c$ and that  $X$ is
arithmetically Cohen-Macaulay if and only if $s = c$. We would like to have
exact sequences of length $c$ also in the case where $X$ is not arithmetically
 Cohen-Macaulay. Thus we have to replace the direct sums of line bundles
 $\cF_i$ by more general sheaves. But we also want to have a uniqueness
 property like for minimal free resolutions. This motivates the following
 concept.

\begin{definition} An exact sequence of sheaves on $\PP^n$
$$
0 \to \cE_s \to \cE_{s-1} \to \ldots \to \cE_1 \to \cJ_X \to 0
$$
is said to be a {\it locally free resolution} of the subscheme $X \subset
\PP^n$ if
\begin{itemize}
\item[(i)]  the sheaves $\cE_i, 1 \leq i
\leq s,$ are locally free (and possibly satisfy some extra conditions),
\item[(ii)] the exact sequence obtained   after cancelation of possibly
occurring redundant line bundles is unique (up to
isomorphisms of exact sequences) among the sequences allowed by (i).
\end{itemize}
Then the latter sequence is  called a {\it minimal locally free resolution}
of $X$.
\end{definition}

Again, we prefer to have locally free resolutions of length $s = c$. Now
consider, for example, that the extra condition in (i) is ``being direct
sums of line bundles''. Then we get back the concept of a (minimal) free
resolution. In this case $s=c$ iff $X$ is arithmetically
 Cohen-Macaulay. Moreover, $s=c$ and $\rank \cE_c = 1$ iff $X$ is
 arithmetically Gorenstein, and $s=c$ and $\rank \cE_1 = c$ iff $X$ is a
 complete intersection.
 The purpose of this paper is to
extend this list of  characterizations of certain arithmetically
Cohen-Macaulay subschemes to certain (locally)
 Cohen-Macaulay subschemes using  locally free resolutions.

Recall that $X$ is arithmetically Cohen-Macaulay if and only if the
intermediate cohomology modules $H^i_*(\cJ_X), i = 1,\ldots,n-c$,
vanish. Thus, the next simplest case from a cohomological point of view
occurs if all the intermediate cohomology modules are annihilated by
the irrelevant maximal ideal $(x_0,\ldots,x_n)$. Then $X$ is called
quasi-Buchsbaum. An even stronger condition is
the property of being arithmetically Buchsbaum. Indeed, $X$ is
arithmetically Buchsbaum if  and only if $X$ and all its consecutive,
sufficiently general hyperplane sections are quasi-Buchsbaum. This is an
important concept which has its origin in a negative answer of Vogel
\cite{Vogel} to a problem of Buchsbaum. It has been stimulating intensive
research (cf., for example, \cite{SV_Amer-J-Math} and in particular
\cite{SV2}). A new
 characterization of arithmetically Buchsbaum subschemes is given by the
 following result (cf.\ also Corollary \ref{geom-char-of-aBM-schemes})
  where $\cOP^i$ denotes the $i$-th exteriour power of the cotangent
  bundle of $\PP^n$. 

\begin{theorem} \label{intro_char_aBM} Let $X \subset \PP^n$ be
  a subscheme of codimension $c$. Then  $X$ is arithmetically Buchsbaum if
  and only if $X$ admits a locally free resolution of the form
$$
0 \to \cF_c \to \ldots \to \cF_2 \to \cF_1 \oplus \bigoplus_j
(\cOP^{p_j}(-e_j))^{s_j} \to \cJ_X \to 0
$$
where $\cF_1, \ldots , \cF_c$ are direct sums of line bundles and $1 \leq
p_j \leq n-c, s_j \geq 0$.
\end{theorem}

The  locally free resolution
in the statement above is called $\Omega$-resolution of $X$.
The theorem extends the main result of Chang in \cite{Chang-Diff-Geom}
from codimension two to arbitrary codimension.  It has been proved as
Corollary II.3.3 in \cite{hab} and independently
by  C.\ Walter
(unpublished) using  different means.
Note that in our approach  the theorem above is a special
case of a more general result characterizing so-called surjective-Buchsbaum
subschemes with finite projective dimension as subscheme of an
arithmetically Gorenstein scheme (cf.\ Theorem \ref{surj-Buchsb-ideals}).
Moreover, we derive necessary conditions for
minimal $\Omega$-resolutions and give some applications.

We can also characterize quasi-Buchsbaum schemes by means of a locally free
resolution provided their dimension does not exceed the codimension.

\begin{theorem} \label{intro_char_quasi-BM}   Let $X \subset \PP^n$ be
  a subscheme of codimension $c \geq \frac{n}{2}$. Then  $X$ is
  quasi-Buchsbaum if
  and only  if $X$ admits a locally free resolution of the form
\begin{eqnarray*}
\lefteqn{
0 \to \cF_c \to \ldots \to \cF_{n-c+1} \to \cF_{n-c} \oplus \bigoplus_j
(\cOP^{2(n-c) -
  1}(-e_{n-c.j}))^{s_{c-n.j}} \to \ldots } & &  \\
& & \to \cF_2 \oplus \bigoplus_j (\cOP^{3}(-e_{2.j}))^{s_{2.j}} \to
 \cF_1 \oplus \bigoplus_j (\cOP^{1}(-e_{1.j}))^{s_{1.j}} \to \cJ_X \to 0
\end{eqnarray*}
where $\cF_1, \ldots , \cF_c$ are direct sums of line bundles and $s_{i.j}
\geq 0$.
\end{theorem}

Again, we give a name to the  locally free resolutions occurring in the
statement. We call them weak $\Omega$-resolutions. The last result
is a special case of Theorem \ref{resolution-by-Omegas}.

Note that the numbers $e_{i.j}$ and $s_{i.j}$ in the  weak
$\Omega$-resolution above  are
uniquely determined by the cohomology of $X$ no matter whether the
resolution is
a minimal one.  Thus, in the case where $X$ is a curve we see that every weak
$\Omega$-resolution is an $\Omega$-resolution. This reflects the well-known
fact that a curve  is
arithmetically Buchsbaum if and only if it is quasi-Buchsbaum.  The
corresponding  statement
fails in higher dimension. For surfaces in $\PP^4$ we will characterize the
weak $\Omega$-resolution of an arithmetically Buchsbaum subscheme. This
result can be applied to conclude from a given weak $\Omega$-resolution
that a quasi-Buchsbaum surface is not arithmetically Buchsbaum.
For example, consider the smooth rational surface $X \subset
\PP^4$ of degree $10$ which can be constructed as degeneracy locus giving
rise  to  the exact sequence
$$
0 \to (\cOP^3(-1))^2 \to \cO(-4)  \oplus (\cOP^1(-3))^2 \to \cJ_X \to 0
$$
(cf.\ \cite{DES-surfaces}, Example B1.15). The second theorem above shows
that $X$ is
quasi-Buchsbaum. However, our results imply that  any surface with such a
resolution cannot be arithmetically Buchsbaum (cf.\ Example
\ref{example-B1.15}).

The paper is organized as follows. In Section 2 we review some results on
duality,
$k$-syzygies and sheaves having at most one non-vanishing intermediate
cohomology. The latter are called Eilenberg-MacLane sheaves. Section 3
describes
$q$-presentations. They are the main tool needed to obtain the locally free
resolutions above.  They have been introduced in local algebra by
Auslander and Bridger \cite{Auslander-Bridger}. We adapt the notion to our
purposes and establish the properties we need later on.

In Section 4  we consider graded modules over a graded Gorenstein
$K$-algebra. First we compare the concepts of a surjective-Buchsbaum,
Buchsbaum and quasi-Buchsbaum module. Second we characterize
surjective-Buchsbaum modules of finite projective dimension by the
existence of a certain locally free resolution.

In Section 5 we consider arithmetically Buchsbaum subschemes of projective
space. Their characterization (Theorem \ref{intro_char_aBM}) is proved
there.  Then we describe how  a free resolution can be obtained from an
$\Omega$-resolution, the behaviour of $\Omega$-resolutions
 under hyperplane sections and consequences for embeddings of abelian
 varieties.  Finally, we
derive restrictions for the degree shifts of the  vector bundles occurring
in a minimal $\Omega$-resolution. This result
 has as an immediate consequence one of the main results of
\cite{Hoa-M}. It says that
the Castelnuovo-Mumford regularity of an arithmetically Buchsbaum subscheme
is almost determined by its index of speciality. As another application our
description of minimal $\Omega$-resolutions has been used  in
\cite{N-gorliaison} to show that new
phenomena occur in the liaison theory of  subschemes of codimension
$\geq 3$ which do not exist in the case where subschemes of codimension $2$
are linked.

In Section 6 we first show that every subscheme  of $\PP^n$ gives rise
to  an exact sequence
$$
0 \to \cE_s \to \ldots \to \cE_1 \to \cJ_X \to 0 \leqno(*)
$$
of length $s \leq \max \{c, \frac{n}{2} \}$
where the sheaves $\cE_i$ are Eilenberg-MacLane sheaves reflecting
the cohomology of $X$ (cf.\ Theorem \ref{resolution_by_E-ML-modules}). It
follows that every equidimensional Cohen-Macaulay subscheme admits a
locally free resolution by Eilenberg-MacLane bundles. From this result we
derive our
characterization of quasi-Buchs\-baum schemes.

In the final section we consider surfaces in $\PP^4$. We interpret
sequence $(*)$ above as saying that every equidimensional Cohen-Macaulay
surface of
codimension two is the degeneracy locus of a morphism between
Eilenberg-MacLane bundles.  We end by establishing a criteria on the  weak
$\Omega$-resolution of a quasi-Buchsbaum surface $X$ which allows to decide
if $X$ is even arithmetically Buchsbaum.

\section{Preliminaries} \label{section_preliminaries}

In this section we fix notation and collect some results which we will need later
on.

Throughout this paper  $R = \oplus_{i \in \mathbb{N}} R_i$  will
denote a graded Gorenstein $K$-algebra where
$R_0$ is the field $K$ and $R$ is generated by the elements
of $R_1$. The
irrelevant maximal ideal $ \oplus_{i > 0} R_i$ of $R$ is denoted by
${\mathfrak m}_R$ or simply ${\mathfrak m}$. Hence $(R,\fm)$ is $^*$local
in the sense of \cite{Bruns-Herzog}.

If $M$ is a module over the graded ring $R$ it is always assumed to be
$\mathbb{Z}$-graded. The set of its homogeneous elements of degree $i$ is
denoted by $M_i$ or $[M]_i$. All homomorphisms between graded $R$-modules
will be morphisms in the category of graded $R$-modules, i.e., will be
graded of degree zero.

If $M$ is an $R$-module $\dim M$ denotes the Krull dimension of $M$. The
symbols $\rank_R$ or simply $\rank$ are reserved to denote the rank of $M$
in case it has one. For a $K$-module $\rankk$ just denotes
 the vector space dimension over the field $K$.

There are two types of duals of an $R$-module $M$ we are going to use. The
$R$-dual of $M$ is $M^* =  \Hom_R(M,R)$
and the $K$-dual $M^{\vee} = \oplus_j \Hom_K ([M]_j, K)$ where $K$ is considered as a
graded module concentrated in degree zero. Both dual modules are
graded $R$-modules.
Note that $R^{\vee}$ is the injective hull of $K^{\vee} \cong
K \cong R/{\mathfrak m}$ in the category of graded $R$-modules. If $\rankk
[M]_i < \infty$ for all integers $i$ then there is a canonical isomorphism
$M \cong M^{\vee \vee}$.
\medskip

Later on we will apply algebraic results in a geometric context.
Thus we mention briefly some relations between the algebraic and geometric
notions which allow us to switch between the two languages.

Let  $\cF$ be a sheaf on
$Z = \Proj (R)$. The cohomology modules of $\cF$ are
$H^i_*({\cF})=\bigoplus_{t\in
  \ZZ}H^i(Z,{\cF}(t))$.

There are two functors relating graded $R$-modules and sheaves of modules
over $Z$. One is the ``sheafification'' functor which associates to each
graded $R$-module $M$ the sheaf $\tilde{M}$. This functor is exact.

In the opposite direction there is the ``twisted global sections'' functor
which associates to each sheaf $\cF$ of modules over $Z$ the graded
$R$-module $H^0_*(\cF)$. This functor is only left exact. If $\cF$ is
quasi-coherent then the sheaf $\widetilde{H^0_*(\cF)}$ is canonically
isomorphic to $\cF$. However, if $M$ is a graded $R$-module then the module
$H^0_*(\tilde{M})$ is not isomorphic to $M$ in general. In fact,
$H^0_*(\tilde{M})$ even
needs not  be finitely generated if $M$ is finitely generated.
However, there is the following comparison result (cf.\ \cite{SV2}).

\begin{proposition} \label{comp} Let $M$ be a graded $R$-module. Then
  there is an exact sequence
$$
0 \to \HH^0(M) \to M \to H^0_*(\tilde{M}) \to \HH^1(M) \to 0
$$
and for all $i \geq 1$ there are isomorphisms
$$
H^i_*(\tilde{M}) \cong \HH^{i+1}(M).
$$
\end{proposition}

It follows, for example, that a coherent sheaf  $\cE$ on $\PP^n$ is locally
free if and only if  all modules $H^i_*(\cE), 1 \leq i < n,$ have   finite
length. As usual we will use ``vector bundle'' and ``locally free sheaf''
interchangeably. Notice that $H^0_*(\tilde{M})$ is isomorphic to
$\displaystyle H^0(M) = \lim_{\stackrel{\longrightarrow}{{\scriptstyle n}}}
\Hom_R(\fm^n,M)$. 
\medskip

Now, we turn to some duality results.  Over the Gorenstein ring $R$ duality
theory is particularly
simple. We denote the index of regularity of a graded ring $A$ by $r(A)$. If
$A$ is just the polynomial ring $K[x_0, \dots , x_n]$ then $r(A) = -n$.  We
will often use the following duality result (cf.\ \cite{S}, \cite{SV2})
without further mentioning.

\begin{proposition} \label{duality} Let $R$ be a graded Gorenstein ring of
  dimension $n+1$. Let $M$ be a graded $A$-module where $A$ is a quotient
  of $R$. Then we have  for all $i \in \mathbb{Z}$ natural isomorphisms of
  graded $R$-modules
$$
H^i_{{\mathfrak m}_A} (M)^{\vee} \cong \Ext^{n+1-i}_R(M,R)(r(R)-1).
$$
\end{proposition}

Let $M$ be an $R$-module of dimension $d$. Then
$$
K_M = \Ext^{n+1-d}_R(M,R) (r(R)-1)
$$
is said to be the {\itshape canonical module} of $M$. It is the module
representing the
functor $\HH^d(M \otimes_R \_\_{})^{\vee}$ (cf.\
\cite{S}). Duality theory also relates the cohomology modules of $M$ and
$K_M$. Later on we will need the following result of Schenzel \cite{S},
Korollar 3.1.3.

\begin{proposition} \label{serredual} Let $M$ be a graded module over the
  Gorenstein ring $R$. Suppose that $\HH^i(M)$ has finite length if $i \neq
  d = \dim M$. Then there are canonical isomorphisms for $i = 2, \ldots ,
  d-1$
$$
\HH^{d+1-i}(K_M) \cong \HH^i(M)^{\vee}.
$$
\end{proposition}

If the assumption on the cohomology of $M$ in the theorem above is
satisfied then $M$ is said to have {\itshape cohomology of finite length}.
A graded  $R$-module has
cohomology of finite length if and only if it is equidimensional and
locally Cohen-Macaulay.
\medskip

Next, we recall some results on  $k$-syzygies. If $M$ is an
$R$-module then $\dim M \leq \dim R$. $M$ is said to be {\itshape maximal}
if $\dim M = \dim R$. Let $Q$ be the total ring of fractions of $R$. Then
an $R$-module $M$ is said to be {\itshape torsion-free} if the natural map
$M \to M \otimes Q$ is injective. The bilinear map $M \times M^* \to R,
(m,\ffi) \mapsto \ffi(m)$, induces a natural homomorphism $h: M \to
M^{**}$. The module $M$ is said to be {\itshape torsionless} if $h$ is
injective, and $M$ is said to be {\itshape reflexive} if $h$ is an
isomorphism. We want to
compare these notions with the following one.

\begin{definition} \label{k-syz} Let $M$ be a graded $R$-module. Then
  $N \neq 0$
  is said to be a {\itshape $k$-syzygy} of $M$ if there is an exact
  sequence of graded $R$-modules
$$
0 \to N \to F_k \stackrel{\ffi_k}{\longrightarrow} F_{k-1} \to  \ldots \to
F_1 \stackrel{\ffi_1}{\longrightarrow} M \to 0
$$
where the modules $F_i, i = 1,\ldots, k,$ are free $R$-modules.

The $R$-module $N$ is
simply called a $k$-syzygy if it is a  $k$-syzygy of some $R$-module.
\end{definition}

Note that a $(k+1)$-syzygy is also a $k$-syzygy.

Let $\ffi: F \to M$ be a homomorphism of $R$-modules where $F$ is
free. Then $\ffi$ is said to be a {\itshape minimal}  homomorphism if $\ffi
\otimes id_{R/\fm} : F/\fm F \to M/\fm M$ is the zero map in case $M$
is free and an
isomorphism in case $\ffi$ is surjective.

In the situation of the definition above $N$ is said to be a {\itshape
  minimal $k$-syzygy} of $M$ if the morphisms $\ffi_i, i = 1, \ldots , k$,
are minimal. It is uniquely determined up to isomorphism by Nakayama's
lemma. If the minimal $k$-syzygy $N$ happens to be free then the exact
sequence in the definition above is called a minimal free resolution of
$M$.

It follows by \cite{Auslander-Bridger} that for a finitely
generated module over a Gorenstein ring the conditions torsionless,
torsion-free and $1$-syzygy are all equivalent. The same applies to the
conditions reflexive  and $2$-syzygy.
It is more difficult to identify third and higher syzygies. For
this we consider the cohomological annihilators $\fa_i(M) =
\Ann_R \HH^i(M)$. Since $R$ is Gorenstein  we have $\dim
R/\fa_i(M) \leq i$ for all integers $i$ where we put $\dim M = - \infty$ if
$M = 0$. Higher
syzygies can be characterized as follows.

\begin{proposition} \label{k-syzch} Let  $M$ be a finitely
generated $R$-module. Then the following conditions are equivalent:
\begin{itemize}
\item[(a)] $M$ is a $k$-syzygy.
\item[(b)] $\dim R/\fa_i(M) \leq i - k$ for all $i < \dim R$.
\end{itemize}
Moreover, if $k \geq 3$ then conditions $\mathrm{(a)}$ and $\mathrm{(b)}$
are equivalent
to the condition that $M$ is reflexive and $\ER^i(M^*,R) = 0$ if $1 \leq i
\leq k-2$.
\end{proposition}

\begin{proof} By local duality we know that the annihilators of
  $\ER^i(M,R)$ and 
  $\HH^{\dim R - i}(M)$ coincide. Hence we get 
$$grade \ER^i(M,R) =
  \dim R - \dim
  R/\fa_{\dim R - i}(M)
$$ 
and the result follows by \cite{Auslander-Bridger},
  Theorem 4.25 and \cite{EG_buch}, Theorem 3.8.
\end{proof}
\medskip

Following Horrocks \cite{Ho2} a maximal  $R$-module $E$ is said to be an
{\itshape Eilenberg-MacLane module} of depth $t$, $0 \leq t \leq n+1 = \dim
R$, if
$$
\HH^j(E) = 0 \quad \mbox{for all} \;  j \neq t \; \mbox{where} \; 0 \leq j
\leq n.
$$

An Eilenberg-MacLane module of depth $n+1$ is Cohen-Macaulay, thus a free
module if it has finite projective dimension. More generally, a relation
between Eilenberg-MacLane modules and syzygy modules is described in the
next result which is proved as Proposition I.3.1 and Theorem I.3.9 in
\cite{hab}.

\begin{proposition} \label{elc} Let $E$ be a
   module of depth $t \leq n$. Then we have:
\begin{itemize}
\item[(a)] If $E$ is an Eilenberg-MacLane
  module with finite projective dimension then $E^*$ is a $(n+2-t)$-syzygy
  of $\HH^t(E)^{\vee}(1-r(R))$.
\item[(b)] If $E$ is reflexive then $E$ is an Eilenberg-MacLane
  module with finite projective dimension if and only if $E^*$ is an
  $(n+2-t)$-syzygy of a module $M$ of dimension $\leq t-2$. In this
  case we have
$$
M \cong \HH^t(E)^{\vee}(1-r(R)).
$$
\end{itemize}
\end{proposition}

A sheaf $\cE$ on $Z = \Proj R$ is said to be an {\it Eilenberg-MacLane
  sheaf \  of depth} $t$ ($1 \leq t \leq n$) if $H^t_*(\cE) \neq 0$
  and $H^i_*(\cE) = 0$ if $i \neq 0, t, n$. Note that in this case
 $H^0_*(\cE)$ is an Eilenberg-MacLane module of depth $t+1$. If $\cE$ is in
addition locally free we call it an Eilenberg-MacLane bundle.

\section{$q$-presentations} \label{section-q-presentations}

In this section we consider $q$-presentations ($q \in \ZZ$)
and describe some of their properties. They will play a crucial
role later on.

These $q$-presentations were first
considered by Auslander and Bridger \cite{Auslander-Bridger} in local
algebra. Their
uniqueness properties were established by  Evans and Griffith
\cite{EG_buch}, \cite{EG_TAMS}. The construction in order to show the
existence of $1$-presentations is due to Serre \cite{Serre-1-presentation}
(cf.\ also \cite{Murthy}). These papers work over a local ring.
However, we stress the fact that the constructions are also
possible over a graded ring.

We slightly modify the notion of $q$-presentations according to our
purposes. We focus on local cohomology and state some useful new properties.

\begin{definition} \label{q_presentation} Let $q$ be an integer with $1
  \leq q \leq \dim R$. An exact sequence of finitely generated, graded
  $R$-modules
$$
0 \to P \stackrel{\ffi}{\longrightarrow} E \to M \to 0
$$
is said to be a {\it $q$-presentation} of  $M$ if
\begin{itemize}
\item[(i)] $P$ has  projective dimension $<q$ \quad and
\item[(ii)] $\HH^j(E) = 0$ for all $j$ with $\dim R - q \leq j < \dim R$.
\end{itemize}
 It is said to be a {\it minimal} $q$-presentation if there does not exist
 a non-trivial free $R$-module $F$ such that $F$ is a direct summand of $P$
 and $E$ and $\ffi$ induces an isomorphism of $F$ onto $F$.
\end{definition}

If the $q$-presentation is not minimal we say that $P$ and $E$ have a
common direct free  summand.

We begin by recalling that a $q$-presentation ``distributes'' the local
cohomology modules of $M$ among $P$ and $E$ (cf.\ \cite{N-gorliaison}).

\begin{lemma} \label{distribute_cohomology} If $0 \to P \to E \to M \to 0$
  is a $q$-presentation of $M$ then
$$
\HH^j(E) \cong \left \{ \begin{array}{ll}
\HH^j(M) & \mif j < \dim R - q \\
0 & \mif \dim R - q \leq j < \dim R
\end{array} \right.
$$
and
$$
\HH^j(P) \cong \left \{ \begin{array}{ll}
0 & \mif j \leq \dim R - q \\
\HH^{j-1}(M) & \mif \dim R - q < j < \dim R.
\end{array} \right.
$$
\end{lemma}

Note that in the definition of a $q$-presentation we have not assumed that
$M$ is maximal nor that $P$ is non-trivial. If $P = 0$ we say that the
$q$-presentation is trivial. Sometimes we know a priori that the modules in
the $q$-presentation of $M$ besides $M$ must be maximal.

\begin{lemma} \label{q-present-has-max-modules} If $1 \leq q < \dim R -
  \dim M$ then $M$ admits a trivial $q$-presentation.

If  $0 \to P \to E \to M \to
0$ is a $q$-presentation then $P$ and $E$ are maximal $R$-modules provided
that
\begin{itemize}
\item[(i)] $1 \leq \dim R - \dim M \leq q \leq \dim R$ \quad or
\item[(ii)] $M$ is maximal and $P \neq 0$.
\end{itemize}
\end{lemma}

\begin{proof} The first claim is immediate by the definition.

Now we show the second one. Assume $\dim R - q \leq \dim M <
\dim R$. Since $\HH^{\dim M}(M) \neq 0$ but $\HH^{\dim M}(E) = 0$ we see
that $P$ is non-trivial. If
$\dim M = \dim R -1$ the exact sequence
$$
0 \to \HH^{\dim R - 1}(M) \to \HH^{\dim R} (P)
$$
shows that $P$ must be maximal. If $\dim M \leq \dim R - 2$ then Lemma
\ref{distribute_cohomology} gives
$$
\HH^i(P) \cong \left \{ \begin{array}{ll}
\HH^{\dim M}(M) \neq 0 & \mif i = \dim M + 1 \\
0 & \mif \dim M + 2 \leq i < \dim R.
\end{array} \right.
$$
Suppose $P$ is not maximal. Then we get $\dim P = \dim M + 1$ and that
$\HH^{\dim P}(P)^{\vee} \cong \HH^{\dim M}(M)^{\vee}$ has dimension $\dim M
< \dim P$ contradicting the fact that the dimension of the canonical module
of $P$ equals $\dim P$. Thus, the second claim follows in case (i).

If $M$ is maximal and $P \neq 0$ then $P$ must be non-trivial by arguing
similarly as above.
\end{proof}

Now we turn to the existence of minimal $q$-presentations.

Since $M$ is by assumption finitely generated it is the epimorphic image of
a finitely generated free $R$-module. This provides a $(\dim R)$-presentation
of $M$ if $M$ has finite projective dimension. In the general
case, a $q$-presentation is constructed by induction on $q$.

Consider the following diagram
$$
\begin{array}{rcl}
& 0 & \\
& \downarrow & \\
& Q & \\
& \downarrow & \\
& G & \\
& \downarrow & \\
0 \to & L & \to \oplus_{j=1}^m R(-e_j) \to M \to 0 \\
& \downarrow & \\
& 0 &
\end{array}
$$
 where the row is a minimal free  presentation of $M$ and the column is a
 minimal
$(q-1)$-presentation of $L$ provided $q \geq 2$. If $q = 1$ we put
$Q = 0$ and $G = L$. In any case, by composition we get a map $G \to
\oplus_{j=1}^m R(-e_j)$ given by elements $g_1, \ldots ,g_m \in G^*$. Now
we choose elements $g_{m+1}, \ldots , g_s \in G^*$ such that $\{g_{m+1},
\ldots , g_s\}$ is a minimal basis of $G^*/(\sum_{j=1}^m g_j R(e_j))$.
Since $L$ is a $1$-syzygy we see that $G$ is a $1$-syzygy too according to
Proposition \ref{k-syzch}. Hence the homomorphism
 $G \to \oplus_{j=1}^{s} R(-e_j)$
 given by $g_{1}, \ldots , g_s$ is injective.  Therefore, using the Snake lemma
we get an exact commutative diagram
$$
\begin{array}{ccccccc}
& 0 & &  0 & &  0 &\\
& \downarrow & & \downarrow & & \downarrow & \\
0 \to & Q  & \to & \oplus_{j=m+1}^{s} R(-e_j) & \to & P & \to 0 \\
& \downarrow & & \downarrow & & \downarrow & \\
0 \to & G & \to & \oplus_{j=1}^{s} R(-e_j) & \to & E & \to 0 \\
& \downarrow & & \downarrow & & \downarrow & \\
0 \to & L & \to & \oplus_{j=1}^{m} R(-e_j) & \to & M & \to 0 \\
& \downarrow & & \downarrow & & \downarrow & \\
& 0 & &  0 & &  0. &\\
\end{array}
$$
A careful analysis shows that its right-hand column is a minimal
$q$-presentation of $M$. This leads to the following result whose
detailed proof can be found in \cite{hab}, Section II.1 (cf.\ also \cite{EG_buch}
 in the local case).

\begin{theorem} \label{q-presentation_exists}  A  finitely generated $R$-module
  $M$ admits for every $q$ with $1 \leq q \leq \dim R$ a minimal $q$-presentation
$$
0 \to P \to E \to M \to 0.
$$
 It is uniquely determined up to isomorphisms of exact sequences.
\end{theorem}

\begin{remark} \label{rem-minimality}
(i) The $1$-presentation of $M$ is minimal if and only if the rank of $P$
equals the number of minimal generators of $\ER^1(M,R)$.

(ii) Contrary to the case $q = 1$ it is not clear how one can decide if a given
$q$-presentation is minimal if $q \geq 2$. A notable exception is described
in Lemma \ref{(c-1)-presentation-of-ideal} later on.
\end{remark}

Finally, we would like to discuss the behaviour of $q$-presentations under
hyperplane sections. Let $l \in R$ be a general linear form and let $0 \to
P \to E \to M \to 0$ be a $q$-presentation of $M$. Suppose $\depth M >
0$. Then the commutative diagram
$$
\begin{array}{ccccccc}
0 \to & P(-1) & \to & E(-1) & \to & M(-1) & \to 0 \\[3pt]
& \downarrow l & & \downarrow l & & \downarrow l &\\[3pt]
0 \to & P & \to & E & \to & M & \to 0. \\
\end{array}
$$
where the vertical maps are  multiplication by $l$ provides the exact
sequence
$$
0 \to \overline{P} \to \overline{E} \to \overline{M} \to 0
$$
of $\overline{R}$-modules where $\overline{R} = R/l R$ and $^-$ denotes the
functor $\_\_\otimes_R \overline{R}$. However, this sequence is not the one
we are looking for. On the one hand it is  in general not a $q$-presentation of
$\overline{M}$ (as $\overline{R}$-module). On the other hand we
would like to have a $q$-presentation of $H^0(\overline{M})$ rather than
$\overline{M}$. This is motivated by applications in geometry. For example,
if $M$ is reflexive then  in general $\overline{M}$ will only be torsion-free
whereas $H^0(\overline{M})$ is a reflexive $\overline{R}$-module. (This
follows by Proposition \ref{k-syzch}). Note also that in case $M$ is an
ideal $I \subset R$ the module $H^0(\overline{I})$ is just the saturation
of $\overline{I}$ in $\overline{R}$. Thus, we will need  the following
technical result later on.

\begin{lemma} \label{q-presentation-and-hyperplane-section} Let
$$
0 \to P \to E \to M \to 0
$$
be a $q$-presentation of $M$ where $1 \leq q \leq \dim
R - 3$. Let
$$
0 \to U \to G \to H^0(\overline{E}) \to 0
$$
be a $q$-presentation of $H^0(\overline{E})$ as $\overline{R}$-module. Suppose
$\depth M > 0$. Then there is an $\overline{R}$-homomorphism $\ffi : G \to
H^0(\overline{M})$ such that
$$
0 \to \ker \ffi \to G \stackrel{\ffi}{\longrightarrow} H^0(\overline{M})
\to 0
$$
is a $q$-presentation of $H^0(\overline{M})$ as $\overline{R}$-module.
\end{lemma}

\begin{proof} As explained above the assumptions provide an exact sequence
$$
0 \to \overline{P} \to \overline{E} \to \overline{M} \to 0
$$
which induces the long exact cohomology sequence
$$
0 \to H^0(\overline{P}) \to H^0(\overline{E}) \to H^0(\overline{M}) \to
H^1(\overline{P}) \to \ldots.
$$
The assumption on $q$ implies $\depth P \geq 4$, thus $\depth \overline{P}
\geq 3$. It follows (cf.\ Proposition
\ref{comp}) that $H^0(\overline{P})
\cong \overline{P}$ and $H^1(\overline{P}) \cong \HH^2(\overline{P}) =
0$. Thus we obtain the following exact commutative diagram where we put $Q = \ker
\ffi$:
$$
\begin{array}{ccccccc}
& & & 0 & & &\\
& & & \downarrow & & &\\
& & & U & & &\\
& & & \downarrow & & &\\
0 \to & Q &  \to & G & \stackrel{\ffi}{\longrightarrow} & H^0(\overline{M}) & \to 0 \\[3pt]
&  & & \downarrow & & \| & \\[3pt]
0 \to & \overline{P} &  \to & H^0(\overline{E}) & \to & H^0(\overline{M}) & \to 0 \\
& & & \downarrow & &  &\\
& & & 0. & & &
\end{array}
$$
The Snake lemma implies the exact sequence
$$
0 \to U \to Q \to \overline{P} \to 0.
$$
Let $F_{\bullet}$ be a minimal free resolution of $P$. Since $l$ is a
non-zero divisor of $R$ and $P$  it is
well-known that $F_{\bullet} \otimes_R \overline{R}$ is a minimal free
resolution of $\overline{P}$ as $\overline{R}$-module. Hence $\overline{P}$
as well as $U$ has projective dimension $<q$ as $\overline{R}$-module. Then the so-called Horseshoe
lemma implies that this is also true for $Q$. Therefore
$$
0 \to Q \to G \to H^0(\overline{M}) \to 0
$$
is a $q$-presentation as claimed.
\end{proof}

\section{Surjective-Buchsbaum modules} \label{section-surj-Buchsbaum}

The theory of Buchsbaum modules started from a negative answer of Vogel
\cite{Vogel} to
a problem posed by  Buchsbaum. The concept has been introduced by
St\"uckrad and Vogel.  We refer to
their monograph \cite{SV2} for a comprehensive
introduction to the subject.
St\"uckrad and Vogel established the so-called surjectivity criterion as a
sufficient condition for a module being Buchsbaum
(cf.\ \cite{SV_Amer-J-Math}, Theorem 1). It
gave rise to the notion of a surjective-Buchsbaum module introduced by
Yamagishi \cite{Yamagishi_surjective-Buchsbaum}.

In this section we give the definition and discuss the preceding
notions. The main result is a characterization of a surjective-Buchsbaum
module of finite projective dimension over a Gorenstein ring with the help
of its $c$-presentation where $c = \dim R - \dim M$.

As before, we restrict our considerations to the graded situation. We
just mention that
all the results have analogues for modules over a local ring containing a
field. The transfer to this situation is obvious and we omit it.

We begin by recalling  some facts of homological algebra. The
inclusion $\soc M
\hookrightarrow \HH^0(M) = \bigcup_{j \geq 1} (0 :_M \fm^j)$ induces natural
homomorphisms of derived functors
$$
\ffi^i_M : \ER^i(K,M) \to \HH^i(M).
$$
Next we need to consider Koszul complexes. For notation and basic properties of
them we refer to \cite{Bruns-Herzog}. Let $\xx = \{ x_1, \ldots , x_s \}$
be a minimal basis of the ideal $\fm$. Then the Koszul complexes
$K_{\bullet}(\xx;M)$ and $K^{\bullet}(\xx;M)$ are up to isomorphism
independent of the chosen minimal basis of $\fm$. Thus it makes sense to
denote them by $K_{\bullet}(\fm;M)$ and $K^{\bullet}(\fm;M)$, respectively,
and the homology modules by $H_i(\fm;M)$ and $H^i(\fm;M)$. Since
$H_0(\fm;R) = R/\fm \cong K$ this isomorphism lifts to a morphism of
complexes from $K_{\bullet}(\fm;R)$ to a minimal free resolution of $K$. It
induces natural homomorphisms
$$
\lambda^i_M : \ER^i(K,M) \to H^i(\fm;M).
$$
Note that $H^0(\fm;M) \cong 0 :_M \fm$ can also be embedded into
$\HH^0(M)$. The induced natural homomorphisms of derived functors are
denoted by
$$
\psi^i_M : H^i(\fm;M) \to \HH^i(M).
$$
Summing up, we have the following commutative diagram for all integers $i$
$$
\xy\xymatrixrowsep{0.5pc}\xymatrix{
\ER^i(K,M) \ar @{->}[dd]^-{\lambda^i_M}  \ar @{->}[dr]^-{\ffi^i_M}  \\
 & \HH^i(M). \\
H^i(\fm;M) \ar @{->}[ur]_-{\psi^i_M} }
\endxy
$$
Now we are ready for the following.

\begin{definition} \label{surjective-Buchsb} The module $M$ is said to be
  {\it surjective-Buchsbaum} if $\ffi^i_M$ is surjective for all $i \neq
  \dim M$. It is said to be {\it Buchsbaum} if $\psi^i_M$ is surjective for
  all $i \neq \dim M$. It is said to be {\it quasi-Buchsbaum} if
$$
\fm \cdot \HH^i(M) = 0 \quad \mbox{for all} ~ i \neq \dim M.
$$
\end{definition}

It is immediate from the commutative diagram above that a
surjective-Buchs\-baum module is Buchsbaum and that every Buchsbaum module is
quasi-Buchs\-baum. Note that these implications are strict in
general. However, if $R$ is regular then $K_{\bullet}(\fm;M)$ is a minimal
free resolution of $K$, i.e., $\ER^i(K,M) \cong H^i(\fm;M)$. Hence, if $R$
is regular then an $R$-module is surjective-Buchsbaum if and only if it is
Buchsbaum. Sometimes the three notions  are
equivalent.

\begin{lemma} \label{hinreichende-Bedingung-fuer-surj-Buchs} Let $M$ be an $R$-module of
  depth $t < \dim M$ such that $\HH^i(M) = 0$ if $i \neq t,
\dim M$. Then $M$ is
  surjective-Buchsbaum if and only if
$$
\fm \cdot \HH^t(M) = 0.
$$
\end{lemma}

\begin{proof} The natural map $\ffi^t_M$ induces an isomorphism 
$$ 
\ER^t(K,M)  \cong \Hom_R(K,\HH^t(M)). 
$$ 
Hence $\ffi^t_M$ is surjective if and only if
  $\fm$ annihilates $\HH^t(M)$.
\end{proof}

Note that the last result applies, for example, if $M$ is the coordinate
ring of a projective curve.
\smallskip

So far the discussion in this section  applies to any graded $K$-algebra
$R$. From now on we will make use of our general
assumption that $R$ is Gorenstein. Moreover, for the rest of this
section we assume that $R$ has positive
dimension.  We put $n+1 = \dim R$ and denote the index of regularity by $r
= r(R)$.

Our next goal  is to
derive a description of the maximal surjective-Buchsbaum modules over $R$
with finite projective dimension. Let
$$
F_{\bullet} : \quad  \ldots  \to F_2 \stackrel{\alpha_2}{\longrightarrow} F_1
\stackrel{\alpha_1}{\longrightarrow} F_0
\stackrel{\alpha_0}{\longrightarrow} K \to 0
$$
be a minimal free resolution of the field $K$. For integers $i$ with $0
\leq i \leq n+1$ we define  $G_i = \coker (\alpha_{n+1-i}^*)(r-1)$.  By
local duality we know  that $\ER^i(K,R) = 0$ if $i \neq n+1$. Thus $G_0(1-r)$
has the following minimal free resolution
$$
0 \stackrel{\alpha_0^*}{\longrightarrow} F_0^*
\stackrel{\alpha_1^*}{\longrightarrow}  \ldots \to F_n^*
\stackrel{\alpha_{n+1}^*}{\longrightarrow} F_{n+1}^* \to G_0(1-r) \to 0.
$$
Moreover, the modules $G_i$ have the following properties.

\begin{lemma} \label{properties-of-indecomposable-max-Buchs} Let $0 \leq i
  \leq n+1$. Then we have:
\begin{itemize}
 \item[(a)] $\HH^0(G_0) \cong K$ and $G_0$ is isomorphic to $K$ if and only
   if $R$ is regular. Otherwise $G_0$ is an Eilenberg-MacLane module.

\item[(b)] If $1 \leq i$ then $G_i$ is a minimal $i$-syzygy of
  $G_0$ and an Eilenberg-MacLane
  module of depth $i$ where $\HH^i(G_i) \cong K$ if $i \leq n$.
\item[(c)] If $i \leq n$ then $G_i^*(r-1)$ is a minimal (n+2-i)-syzygy of
  $K$ whereas $G_{n+1}^*(r-1) = R$.
\item[(d)] $G_i$ is surjective-Buchsbaum, indecomposable and has finite
  projective dimension.
\end{itemize}
\end{lemma}

\begin{proof} According  to the construction of $G_0$ we have $\ER^i(G_0,R) = 0$
  if $1 \leq i \leq n$ and $\ER^{n+1}(G_0,R) \cong K(r-1)$. Note also that
  $K$ has finite projective dimension if and only if $R$ is regular. Thus
  claim (a) follows by local duality.

Claims (b) and (c) are immediate from the resolutions of $K$ and $G_0$
above (cf.\ also Proposition \ref{elc}).

Suppose $G_0$ is decomposable and maximal. Since $K$ is indecomposable it
follows by (a) that
one of the direct summands of $G_0$ must be maximal
Cohen-Macaulay, thus free. Hence $G_0^*$ has a free direct summand
contradicting the fact that $G^*_0$ is a minimal syzygy of the indecomposable
module $K(1-r)$. The remaining assertions of (d) follow using Lemma
\ref{hinreichende-Bedingung-fuer-surj-Buchs}.
\end{proof}

The next result implies in particular that up to degree shift the modules
$G_i\ (0 \leq i \leq n+1)$ are the only indecomposable maximal surjective-Buchsbaum modules of
finite projective dimension with the exception of $G_0$ in the case where
$R$ is regular.

\begin{proposition} \label{Char-of-max-surj-Buchsbaum} Let $M$ be a maximal
   module with finite projective dimension. Then the following conditions
   are equivalent
\begin{itemize}
\item[(a)] $M$ is surjective-Buchsbaum.
\item[(b)]
$$
M \cong F \oplus \bigoplus_{i=0}^n \bigoplus_{j \in \ZZ} G_i^{s_{ij}}(-j)
$$
where $F$ is free and $s_{ij} = [\hh^i(M)]_j$.
\end{itemize}
\end{proposition}

\begin{proof} Since $R$ is Gorenstein the $R$-module $M$ has finite
  projective dimension if and only it has finite injective dimension. Thus
  the claim follows from the proof of \cite{Kawasaki-surj-Buchs}, Theorem
  3.1.
\end{proof}

Observe that the direct sum in the last statement is certainly finite
because $M$ is finitely generated by assumption.

The quoted result of Kawasaki has its origin in a theorem of Goto
\cite{Goto-max-Buchsbaum} who
proved the statement above in the case where $R$ is regular (cf.\ also
\cite{EG},
Theorem 3.2). It is much more difficult to describe the
isomorphism classes of maximal modules of finite projective dimension which
are Buchsbaum or quasi-Buchsbaum. In particular, there are infinitely many
of them if $R$ is not regular (cf.\ \cite{Yoshino_max-Buchsbaum})

We need one more preparatory result.

\begin{lemma} \label{transfer-in-c-present} Let
$$
0 \to P \to E \to M \to 0
$$
  be a $q$-presentation where $1 \leq q \leq \dim R - \dim M$. Then $M$ is
  surjective-Buchsbaum respectively quasi-Buchsbaum if and only if $E$ has
  the corresponding property.
\end{lemma}

\begin{proof} The given $q$-presentation induces for all integers $i$ the
  following  commutative diagram with exact rows:
\begin{equation*}
\begin{CD}
\ER^i(K,P) @>>> \ER^i(K,E)  @>>> \ER^i(K,M)  @>>> \ER^{i+1}(K,P) \\
@VV{\ffi^i_P}V @VV{\ffi^i_E}V @VV{\ffi^i_M}V @VV{\ffi^{i+1}_P}V \\
\HH^i(P)  @>>>  \HH^i(E)  @>>> \HH^i(M)  @>>> \HH^{i+1}(P). \\
\end{CD}
\end{equation*}
According to the definition of a $q$-presentation we  get $\depth P \geq
\dim R + 1 - q =
n+2-q$. Hence the left- and the right-hand side of the  rows in  the diagram
above vanish if $i \leq n-q$. It follows that for these $i$ the map
$\ffi^i_E$ is surjective if and only if $\ffi^i_M$ is surjective. Since
$\dim M \leq n+1-q$ by assumption and $\HH^i(E) = 0$ if $n+1-q \leq i \leq
n$ by Lemma \ref{distribute_cohomology} we obtain that $E$ is
surjective-Buchsbaum if and only if $M$ is. This also shows the
corresponding assertion with respect to the quasi-Buchsbaum property by
just considering the lower row in the diagram above.
\end{proof}

After these preparations the main result of this section follows
easily. Recall that the canonical module of $M$ is denoted by $K_M$.

\begin{theorem} \label{char-of-surj-Buchsbaum} Let $M$ be an
   $R$-module such that  $c = \dim R - \dim M > 0$.  Then the
  following conditions are equivalent:
\begin{itemize}
\item[(a)] $M$ is surjective-Buchsbaum with finite projective dimension.
\item[(b)] $M$ admits a $c$-presentation
$$
0 \to P \to E \to M \to 0
$$
whereby  $P$ has projective dimension $c-1$, $P^*$ is a $c$-syzygy of $K_M(1-r)$ and
$$
E \cong F \oplus \bigoplus_{i=0}^{\dim M -1} \bigoplus_{j \in \ZZ} G_i^{s_{ij}}(-j)
$$
where $F$ is free and $s_{ij} = [\hh^i(M)]_j$.
\end{itemize}
\end{theorem}

\begin{proof} Let $0 \to P \to E \to M \to 0$ be a $c$-presentation of
  $M$. It follows by Lemma \ref{q-present-has-max-modules} that
  $P$ and $E$ must be maximal modules. Hence Lemma
  \ref{distribute_cohomology} implies that $P$ is an Eilenberg-MacLane
  module of depth $n+2-c$ where $\HH^{n+2-c}(P) \cong \HH^{n+1-c}(M) \neq
  0$ because of $\dim M = n+1-c$. Now the Auslander-Buchsbaum formula gives
  the claim on the projective dimension and Proposition \ref{elc} yields
  that $P^*$ is a $c$-syzygy of $\HH^{n+1-c}(M)^{\vee}(1-r) \cong
  \ER^c(M,R) = K_M(1-r)$.

Since $P$ has finite projective dimension the same is true for $E$
if and only if $M$ has finite projective dimension.

Thus the asserted equivalence follows by Proposition
\ref{Char-of-max-surj-Buchsbaum} and
Lemma \ref{transfer-in-c-present}.
\end{proof}

Observe that the module $P$ in the $c$-presentation above is
torsion-free but not reflexive. Replacing $P$ by its minimal free
resolution we obtain.

\begin{corollary} \label{sur-B_res}  Let $M$ be an
   $R$-module such that  $c = \dim R - \dim M > 0$.  Then the
  following conditions are equivalent:
\begin{itemize}
\item[(a)] $M$ is surjective-Buchsbaum with finite projective dimension.
\item[(b)] $M$ admits a locally free resolution
$$
0 \to F_c \to \ldots \to F_1 \to F_0 \oplus \bigoplus_{i=0}^{\dim M -1}
\bigoplus_{j \in \ZZ} G_i^{s_{ij}}(-j) \to M \to 0
$$
where $F_0,\ldots, F_c$ are free $R$-modules.
\end{itemize}
\end{corollary}

\section{Arithmetically Buchsbaum subschemes}
\label{section-arith-Buchsb-subschemes}

In this section we show how the methods of the previous sections can be applied
in order to study an arithmetically Buchsbaum subscheme of projective
space. We begin by considering the $(c-1)$-presentation of an ideal $I \subset
R$ because it turns out that it is more useful than the $c$-presentation of
$R/I$. Then we obtain a characterization of projective arithmetically
Buchsbaum subschemes by means of the so-called $\Omega$-resolution.

Throughout the rest of this paper $I$ will be a homogeneous ideal in the graded
Gorenstein $K$-algebra $R$. We put $A = R/I, n+1 = \dim R, r = r(R)$ and
denote the
codimension of $I$ by $c$, i.e. $\dim A = n+1-c$. We will always assume
that $2 \leq c \leq  n$. Otherwise the results become rather trivial.

Our first goal is a minimality criterion for the $(c-1)$-presentation of
$I$. In order to unify the statement we say that a free module $F$ is a
{\it $0$-syzygy} of a module $N$ if there is an epimorphism $\ffi : F \to
N$. $F$ is called a minimal $0$-syzygy if the map $\ffi$ is minimal.

\begin{lemma} \label{(c-1)-presentation-of-ideal} Let
$$
0 \to P \to E \to I \to 0
$$
be a $(c-1)$-presentation of $I$. Then $P$ is a reflexive Eilenberg-MacLane
module of projective dimension $c-2$ and $P^*$ is a $(c-1)$-syzygy of
$K_A(1-r)$.

Moreover, the given $(c-1)$-presentation is minimal if and only if $P^*$ is
a minimal $(c-1)$-syzygy of $K_A(1-r)$.
\end{lemma}

\begin{proof} Lemma \ref{distribute_cohomology} provides that $P$ is an
  Eilenberg-MacLane module of depth $n+3-c$ where $\HH^{n+3-c} \cong
  \HH^{n+1-c}(A)$ if $c \geq 3$. This implies that $P$ is reflexive and
  has  the claimed projective dimension. Proposition
  \ref{elc} gives that $P^*$ is
  a a $(c-1)$-syzygy of $K_A(1-r)$ if $c \geq 3$. If $c=2$ the last claim is
  immediate from the given $1$-presentation and the minimality assertion
  follows by Remark \ref{rem-minimality}.

It remains to show the minimality criterion if $c \geq 3$. Assume that the
given presentation is not minimal. Then there is a non-trivial free module
$F$ such that $P \cong F \oplus Q$ for some module $Q$. Since $Q$ has the
same intermediate cohomology as $P$, $Q^*$ is also a $(c-1)$-syzygy of $K_A(1-r)$ by
Proposition \ref{elc}, i.e.\ $P^*$ is not a minimal $(c-1)$-syzygy of
$K_A(1-r)$.

Now suppose that the given $(c-1)$-presentation is minimal. We have to show
that $P^*$ does not have a free direct summand. Dualizing the presentation
we obtain the exact sequence
$$
0 \to R \to E^* \to P^* \to \ER^1(I,R) \to \ldots .
$$
But $\ER^1(I,R) \cong \ER^2(A,R) \cong \HH^{n-1}(A)^{\vee}(1-r) = 0$
because $\dim A = n+1-c \leq n-2$. Thus we see: if $P^*$ has a free direct
summand then $E^*$ has this direct summand too, contradicting the assumed
minimality of the given presentation of $I$.
\end{proof}

Using the modules $G_i$ defined before Lemma
\ref{properties-of-indecomposable-max-Buchs} we state the main result of
this section.

\begin{theorem} \label{surj-Buchsb-ideals}  With the notation above the
  following conditions are equivalent
\begin{itemize}
\item[(a)] $A = R/I$ is surjective-Buchsbaum of finite projective dimension.
\item[(b)] $I$ admits an exact sequence
$$
0 \to F_c \to \ldots \to F_2 \to F_1 \oplus \bigoplus_{i=1}^{n+1-c}
\bigoplus_{j \in \ZZ} (G_i(-j))^{t_{ij}} \to I  \to 0
$$
where $F_1, \ldots , F_c$ are free and $t_{ij} = \rankk [\HH^{i}(I)]_j$.
\item[(c)]  $I$ admits an exact sequence
$$
0 \to F_c' \oplus \bigoplus_{i=c}^{n} \bigoplus_{j
  \in \ZZ} (G_i(-j))^{t_{ij}'}  \to F_{c-1}' \to \ldots \to F_1'  \to I  \to 0
$$
where $F_1', \ldots , F_c'$ are free and $t_{ij}' = \rankk
[\HH^{i+1-c}(I)]_j$.
\end{itemize}
\end{theorem}

\begin{proof} We begin with showing the equivalence of (a) and (b).
Consider a
$(c-1)$-presentation of the ideal $I$:
$$
0 \to P \to E \to I \to 0.
$$
It induces the following commutative diagram with exact rows:
\begin{equation*}
\begin{CD}
\ER^i(K,P) @>>> \ER^i(K,E)  @>>> \ER^i(K,I)  @>>> \ER^{i+1}(K,P) \\
@VV{\ffi^i_P}V @VV{\ffi^i_E}V @VV{\ffi^i_I}V @VV{\ffi^{i+1}_P}V \\
\HH^i(P)  @>>>  \HH^i(E)  @>>> \HH^i(I)  @>>> \HH^{i+1}(P). \\
\end{CD}
\end{equation*}
Since $\depth P = n+3-c$ we see that for $i \leq n+1-c$ the map
$\ffi^i_E$ is surjective if and only if $\ffi^i_I$ is surjective. The same
reasoning starting with the exact sequence
$$
0 \to I \to R \to A \to 0
$$
shows that for $i < n$ the map $\ffi^i_A$ is surjective if and only if
$\ffi^{i+1}_I$ is. It follows that $A$ is surjective-Buchsbaum if and only if
$\ffi^i_I$ is surjective for all $i \leq n+1-c$.
But by the definition of a $(c-1)$-presentation we have for all $i$ with
$n+2-c \leq i \leq n$ that $\HH^i(E) = 0$. This shows that $A$ is
surjective-Buchsbaum if and only if $E$ is.

By the previous lemma $P$ has projective dimension $c-2$.
Thus the assertion on the integers $t_{ij}$ follows by Proposition
\ref{Char-of-max-surj-Buchsbaum} and
Lemma \ref{properties-of-indecomposable-max-Buchs}.

In order to show the equivalence of (a) and (c) we consider the beginning
of a free resolution of $X$:
$$
0 \to N  \to F_{c-1}' \to \ldots \to F_1'  \to I  \to 0.
$$
According to Proposition \ref{Char-of-max-surj-Buchsbaum} we have to show
that $N$ is surjective-Buchsbaum of finite projective dimension if and only
if $I$ has this property. This
follows by shopping the exact sequence into short exact sequences and by
using the following fact: Let
$$
0 \to N \to F \to M \to 0
$$
be an exact sequence of maximal modules where $F$ is free and $\HH^n(M) =
0$. Then we have $\HH^i (N) \cong \HH^i (M)$ for all $i \leq n$ and $N$ is
surjective-Buchsbaum of finite projective dimension if and only if $M$ has
this property. This fact is proved using  similar arguments as in the first
part of the proof. We leave the details to the reader.
\end{proof}

\begin{remark} \label{surj-Buchsb-under-liaison} Let $R$ be a ring which
  is not regular. Let $I \subset R$ be a saturated ideal which is
  surjective-Buchsbaum of finite projective dimension but not perfect. Let
  $\fc \subset I$
  be a Gorenstein ideal with the same codimension as $I$. Then the ideal $J
  = \fc : I$ is linked to $I$ and has infinite projective
  dimension. Indeed, according to out last result and \cite{N-gorliaison},
  Proposition 3.8 the ideal $J$ admits an exact sequence
$$
0 \to F_c' \oplus \bigoplus_{i=1}^{n+1-c}
\bigoplus_{j \in \ZZ} (G_i^*(+j))^{t_{ij}} \to F_{c-1}' \to   \ldots  \to F_1'  \to J(s)  \to 0
$$
where $s$ is a suitable integer and the modules $F_1',\ldots,F_c'$ are
free. The modules $G_i^*$ are syzygy modules of the residue field
$K$. Since $R$ is not regular by assumption, it follows that the $G_i^*$ have
infinite projective dimension. Thus $J$ has this property too due to the
exact sequence above.

Hence we have shown that the property of being surjective-Buchsbaum with
finite projective dimension is not preserved in the whole liaison class of
such an ideal $I$ if $R$ is not regular. On the other hand the property is
preserved in the whole even liaison  of $I$ according to
\cite{N-gorliaison}, Theorem 3.10.
\end{remark}

Now we  want
to formulate the consequences of Theorem \ref{char-of-surj-Buchsbaum}  for
an arithmetically Buchsbaum subscheme $X \subset \PP^n = \PP^n_K$.
For the rest of this section we assume the field
$K$ to be infinite. We denote by $I(X) \subset R$ the homogeneous ideal of
$X$. Then $A = R/I(X)$ is the homogeneous coordinate
ring of $X$. Recall that a subscheme $X$ is said to be
arithmetically Buchsbaum if its homogeneous coordinate ring is
Buchsbaum. As usual we denote the ideal sheaf of $X$ by $\cJ_X$.

Then we have.

\begin{corollary} \label{geom-char-of-aBM-schemes} Let $X \subset \PP^n$ be
  a subscheme of codimension $c$. Then the following conditions are
  equivalent:
\begin{itemize}
\item[(a)] $X$ is arithmetically Buchsbaum.
\item[(b)] $\fm \cdot H^i_*(\cJ_{X \cap L}) = 0$ for any linear subspace $L
  \subset \PP^n$ of dimension $> c$ intersecting $X$ transversally and
  all $i$ with $1 \leq i \leq  \dim L - c$.
\item[(c)] There is an exact sequence
$$
0 \to \cF_c \to \ldots \to \cF_2 \to \cF_1 \oplus \bigoplus_j
(\cOP^{p_j}(-e_j))^{s_j} \to \cJ_X \to 0
$$
where $\cF_1, \ldots , \cF_c$ are direct sums of line bundles and $1 \leq
p_j \leq n-c$.
\item[(d)] There is an exact sequence
$$
0 \to \cF_c'  \oplus \bigoplus_j
(\cOP^{p_j'}(-e_j'))^{s_j'} \to \cF_{c-1}'  \to \ldots  \to \cF_1' \to \cJ_X
\to 0
$$
where $\cF_1', \ldots , \cF_c'$ are direct sums of line bundles and $c \leq
p_j' \leq n-1$.
\end{itemize}
\end{corollary}

\begin{proof} The equivalence of (a) and (b) is shown in \cite{SV2}. Hence
  the claim follows by Theorem \ref{surj-Buchsb-ideals} using $H^i_*(\cJ_X)
  \cong \HH^{i+1}(I(X))$ if $i \geq 1$ and $\tilde{G_i} \cong \cOP^{i-1}$
  by Lemma \ref{properties-of-indecomposable-max-Buchs}.
\end{proof}

\begin{remark} (i)  The last statement is a generalization of the main
  result of Chang in \cite{Chang-Diff-Geom}. She proved it for subschemes of
  codimension $2$. The result has been proved independently by
  C.\ Walter (unpublished).

(ii) Since $(\cOP^p)^* \cong \cOP^{n-p}(n+1)$ the corollary shows in
particular that any arithmetically Buchsbaum subscheme of codimension $c$
is the zero scheme of a section of a vector bundle $\cE = \oplus_j
(\cOP^{q_j}(f_j))^{s_j}$ with $c \leq q_j \leq n$. This section cannot be a
general one if $\cE$ is not a direct sum of line bundles because then $\cE$
has rank $\geq n$.
\end{remark}

Following Chang we want to introduce a name for the exact sequence
described in part (c) of the corollary above.

\begin{definition} \label{Def-Omega-resolution} A subscheme $X$ of
  codimension $c$ is said to have an {\it $\Omega$-resolu\-tion} if there
  exists an exact sequence
$$
0 \to \cF_c \stackrel{\alpha_c}{\longrightarrow} \ldots \to \cF_2 \stackrel{\alpha_2}{\longrightarrow} \cF_1 \oplus \bigoplus_j
(\cOP^{p_j}(-e_j))^{s_j} \to \cJ_X \to 0
$$
where $\cF_1, \ldots , \cF_c$ are finite direct sums of line bundles, $1 \leq
p_j \leq n-c, 1 \leq s_j$ and $(p_j,e_j)$ are all distinct ordered pairs of
integers.

The $\Omega$-resolution is said to be {\it minimal} if there is no line
bundle $\cL$ in the resolution such that the restriction of $\alpha_i$ ($2
\leq i \leq c$) to
$\cL$ induces an isomorphism of $\cL$ onto $\cL$.
\end{definition}

Corollary \ref{geom-char-of-aBM-schemes} says that a subscheme $X$ has an
$\Omega$-resolution if and only if it is arithmetically Buchsbaum.

It is not difficult to see that the numbers $p_j, e_j, s_j$ in the
$\Omega$-resolution are uniquely determined by $X$ because
$h^{p_j}(\cJ_X(e_j)) = s_j$ are the only non-zero intermediate cohomology
groups
of $\cJ_X$. In fact, there is a stronger uniqueness property which shows
that an $\Omega$-resolution is indeed a locally free resolution in the
sense  of  the introduction.

\begin{lemma} \label{Omega-res-is-unique}  The minimal $\Omega$-resolution
  of an arithmetically Buchsbaum subscheme is uniquely determined (up to isomorphism).
\end{lemma}

\begin{proof} Every minimal $\Omega$-resolution of $X$ gives rise to a
  minimal  $(c-1)$-presentation $0 \to P \to E \to I(X) \to 0$  and a
  minimal free resolution of $P$ and vice versa. But a minimal
  $(c-1)$-presentation is
  uniquely determined due to  Theorem \ref{q-presentation_exists}.
\end{proof}

\begin{remark} \label{proof_of_thm1_intro} Putting together Corollary
  \ref{geom-char-of-aBM-schemes} and
  Lemma \ref{Omega-res-is-unique} we obtain Theorem \ref{intro_char_aBM}
  of the introduction.
\end{remark}

Now we want to draw some consequences of the existence of a
$\Omega$-resolution.
Since the Koszul complex provides a  minimal free resolution of $\cOP^p$
the mapping cone construction implies  the
following information on the free resolution of $X$.

\begin{corollary} \label{resolution-derived-from-Omega-res} If $X$ has an
  $\Omega$-resolution as in Definition \ref{Def-Omega-resolution} then
  $\cJ_X$ has a (possibly non-minimal) free resolution of the form:
\begin{eqnarray*}
\lefteqn{0 \to \bigoplus_{p_j = 1} (\cO(-e_j-n-1))^{s_j} \to
  \bigoplus_{1 \leq p_j 
  \leq 2} (\cO(-e_j-p_j-n+1))^{s_j \binom{n+1}{p_j+n-1}} \to \ldots
 } \\ 
& & \to \bigoplus_{1 \leq p_j \leq n-c}
 (\cO(-e_j-p_j-c-1))^{s_j \binom{n+1}{p_j+c+1}} \\
& & \to \cF_c \oplus \bigoplus_j (\cO(-e_j-p_j-c))^{s_j
  \binom{n+1}{p_j+c}}\hspace*{5.0cm}  \\
& & \to \ldots
\to \cF_1 \oplus \bigoplus_j (\cO(-e_j-p_j-1))^{s_j \binom{n+1}{p_j+1}} \to
\cJ_X \to 0.
\end{eqnarray*}
\end{corollary}

\begin{remark} \label{end-of-resolution} Let $t$ be the depth of $A = R/I(X)$. Suppose $t < \dim
  A$. Then a slight generalization of a result of Rao \cite{Rao-Invent},
  Theorem 2.5 implies that the last free module in the minimal free
  resolution of $A$ equals the last free module in the minimal free
  resolution of $\HH^t(A)$ if the latter has finite length. In case $A$ is
  Buchsbaum a lot more is true. The previous corollary implies that the last
  $n+1-t-c$ free modules in the minimal free resolution of $\cJ_X$ are
  completely determined by the intermediate cohomologies of $X$ because at
  these spots  cancellation is impossible. This
  means using the previous notation
$$
\TR_i(K,A) \cong \bigoplus_j (K(-e_j-p_j-i))^{s_j \binom{n+1}{p_j+i}} \quad
\mif i \geq c+1.
$$
\end{remark}

Next, we want to describe the $\Omega$-resolution of the general hyperplane
section of a subscheme having an $\Omega$-resolution.

\begin{lemma} \label{Omeag-res-unter-hyperpl-sect} Let $X$ be a subscheme
  having an $\Omega$-resolution as in Definition
  \ref{Def-Omega-resolution}. Let $H \subset \PP^n$ be a general
  hyperplane. Then $X \cap H$ has an $\Omega$-resolution as follows:
\begin{eqnarray*}
0 \to \cF_c|_H \oplus \bigoplus_{p_j=n-c}(\cO_H(-e_j-n))^{s_j
\binom{n}{0}} \to \ldots \hspace*{4.5cm} \\
\to \cF_2|_H \oplus \bigoplus_{p_j=n-c}(\cO_H(-e_j-n-2+c))^{s_j 
  \binom{n}{c-2}} \\
\to \cF_1|_H \oplus \bigoplus_{p_j=1}(\cO_H(-e_j-1))^{s_j} \oplus
\bigoplus_{p_j=n-c}(\cO_H(-e_j-n-1+c))^{s_j \binom{n}{c-1}} \oplus \\
\bigoplus_{1 \leq p_j < n-c} (\Omega_H^{p_j}(-e_j))^{s_j} \oplus
\bigoplus_{2 \leq p_j \leq n-c} (\Omega_H^{p_j-1}(-e_j-1))^{s_j} \\
\to \cJ_{X \cap H} \to 0.
\end{eqnarray*}
\end{lemma}

\begin{proof} We will follow the approach described in Lemma
  \ref{q-presentation-and-hyperplane-section}. It is well-known that
$$
\cOP^{p}|_H \cong \Omega_H^p \oplus \Omega_H^{p-1}(-1).
$$
Write the second
  term in the given $\Omega$-resolution of $X$ as $\cE \oplus
  \bigoplus_{p_j=n-c} (\cOP^{n-c}(-e_j))^{s_j}$. Then its restriction to
  $H$ can
be obtained by applying the last formula. We have to find a
$(c-1)$-presentation of this restriction. Fortunately this is easy. We
obtain by resolving $\bigoplus_{p_j=n-c} (\Omega_H^{n-c}(-e_j))^{s_j}$ the
$(c-1)$-presentation
\begin{eqnarray*}
\lefteqn{0 \to \bigoplus_{p_j=n-c} (\Omega_H^{n-c+1}(-e_j))^{s_j} \to} \\
& & \cE|_H  \oplus
\bigoplus_{p_j=n-c} (\Omega_H^{n-c-1}(-e_j-1))^{s_j} \oplus
\bigoplus_{p_j=n-c} (\cO_H(-e_j-n-1+c))^{s_j \binom{n}{c-1}} \to \\
& &  \cE|_H \oplus
\bigoplus_{p_j=n-c} (\Omega_H^{n-c-1}(-e_j-1))^{s_j} \oplus
\bigoplus_{p_j=n-c} (\Omega_H^{n-c}(-e_j))^{s_j} \to 0.
\end{eqnarray*}
The Koszul complex provides a resolution of the term on the left-hand
side. Thus an application of the Horseshoe lemma gives the desired
$\Omega$-resolution as in Lemma
\ref{q-presentation-and-hyperplane-section}.
\end{proof}

Roughly speaking, we get the $\Omega$-resolution of the hyperplane section
by restricting any bundle occurring in the $\Omega$-resolution of $X$ to $H$
and then replacing $\bigoplus_{p_j=n-c} (\Omega_H^{n-c}(-e_j))^{s_j}$ by its
free resolution.

Cutting $X$ by $n-c$ general hyperplanes we obtain a zero-dimensional
subscheme $Z$. Repeated use of the lemma above provides an
$\Omega$-resolution of $Z$ which is in fact a free resolution reflecting
the fact that $Z$ is arithmetically Cohen-Macaulay.
\smallskip


Next, we want to consider abelian varieties. We need the following
sufficient but not necessary Buchsbaum criteria (cf., for example,
\cite{SV2}, Proposition I.3.10) which generalizes Lemma
\ref{hinreichende-Bedingung-fuer-surj-Buchs}.

\begin{lemma} \label{suff-crit-fuer-Buchsbaum} Let $X \subset \PP^n$ be a
  subscheme such that the following two conditions are satisfied:
\begin{itemize}
\item[(i)] $\fm \cdot H^i_*(\cJ_X) = 0$ for all $i$ with $1 \leq i \leq
  \dim X$.
\item[(ii)] If $H^i(\cJ_X(e)) \neq 0$ and $H^j(\cJ_X(f)) \neq 0$ for
  some integers $i, j$ with $1 \leq  i < j \leq \dim X$ then $(i+e) - (j+f)
  \neq 1$.
\end{itemize}
Then $X$ is arithmetically Buchsbaum.
\end{lemma}

\begin{example} Let $\cL$ be an ample line bundle on an abelian variety
  $Y$ of dimension $g$. Then $\cL^m$ is very ample if $m \geq 3$ according
  to a result of Lefschetz (cf.\ \cite{Mumford-Tata}, p.\ 163).  In this case
  $\cL^m$ provides an embedding of
  $Y$ into $\PP( H^0(Y,\cL^m))$. Let us denote the image by $X$. Then $X$ is
  projectively normal due to Koizumi and Ohbuchi
  (cf.\ \cite{Lange-Birkenhake}, Theorem 7.3.1) and its
  homogeneous ideal is generated by quadrics if $m \geq 4$ according to
  Kempf \cite{Kempf}.  From the results on the cohomology of line
  bundles on an 
  abelian variety \cite{Mumford-Tata}, Section 16  it follows that
$$
H^i_*(\cJ_X) \cong K^{\binom{g}{i-1}} \quad \mif 2 \leq i \leq g.
$$
Hence if $m \geq 4$ then  $X$ is arithmetically Buchsbaum by the lemma
above and has an 
$\Omega$-resolution of the form:
$$
0 \to \cF_c \to \ldots \to \cF_2 \to \cO_{\PP^n}^{\alpha}(-2) \oplus \bigoplus_{2 \leq
  i \leq g} (\cOP^i)^{\binom{g}{i-1}} \to \cJ_X \to 0
$$
where $c = h^0(Y,\cL^m) - 1 - g$ and $\alpha > 0$. Using the fact that $X$
is subcanonical this resolution can be described a little more
detailed. It follows for example that $\cF_c = \cO (-g-c)$ if the
$\Omega$-resolution above is minimal.

In \cite{Pareschi} Pareschi has obtained information on the beginning of the
minimal free resolution of $X$ by showing that it is linear in the first
stages. The $\Omega$-resolution of $X$
above allows to determine the end of 
the minimal free resolution of $X$ (cf. Remark \ref{end-of-resolution}). It
follows immediately that  the whole free resolution of $X$ cannot be linear.
\end{example}

Finally we want to collect some information on the twists of the line
bundles occurring in a minimal $\Omega$-resolution. We will use the
following notation for a graded  $R$-module $N$:
$$
e(N) = \sup \{j \in \ZZ \s [N]_j \neq 0 \}
$$
and
$$
e^+(N) = e(N/\fm N).
$$
Note that $e^+(N)$ is just the maximal degree of a
minimal generator of $N$ if $N$ is finitely generated.

\begin{proposition} \label{shifts-in-Omega-res} Suppose $X$ has a minimal
  $\Omega$-resolution as in Definition \ref{Def-Omega-resolution}. Let
$$
e(X) = e(H^{n-c+1}_*(\cJ_X))
$$
be the index of speciality and let $\cF_i = \oplus_k \cO_{\PP^n}(-d_{ik}) \
(1 \leq i \leq c)$. Then it holds:
\begin{itemize}
\item[(a)] $$1 \leq p_j + e_j \leq e(X) + n + 2 - c \quad \mbox{for all} \;
  j. $$
\item[(b)] $$\min_{j,k} \{ d_{1k}, p_j + e_j \} - 1 \leq d_{ik} - i \leq
  e(X) + n +
  1 -c \quad \mif 1 \leq i \leq c. $$
\end{itemize}
\end{proposition}

\begin{proof} Taking global sections the given $\Omega$-resolution provides
  a minimal $(c-1)$-presentation of $I = I(X)$:
$$
0 \to P \to E \to I \to 0
$$
and a minimal free resolution of $P$:
$$
0 \to F_c \to \ldots \to F_2 \to P \to 0.
$$
First let us suppose $c \geq 3$. Since $\ER^i(P,R) = 0$ if $i > 0$ and $i
\neq c-2$ dualizing the resolution of $P$ with respect to $R$ gives the
exact sequence:
$$
0 \to P^* \to F_2^* \to \ldots \to F_c^* \to \ER^{c-2}(P,R) \to
0. \leqno(*)
$$
Moreover, we get the exact sequence
$$
0 \to R \to E^* \to P^* \to 0.
$$
Let $y \in E^*$ be the image of the unit element of $R$. Then $I$ is the
order ideal of $y$. Since $E^*$ is a $(c+1)$-syzygy by the explicit
description of $E$ (cf.\ Corollary \ref{geom-char-of-aBM-schemes}) and $I$ has codimension $c$ it follows by
\cite{EG_buch}, Theorem 3.14 that $y$ is not a minimal generator of
$E^*$. This implies in particular that $a(E^*) = a(P^*)$ and $e^+(E^*) =
e^+(P^*)$. Now we have
$$
(\tilde{E})^* \cong \cF_1^* \oplus
\bigoplus_j(\cOP^{n-p_j}(n+1+e_j))^{s_j}.
$$
It follows
$$
a(E^*) = \min_{j,k} \{ -d_{1k}, -p_j - e_j \} \quad \mbox{and} \quad e^+(E^*) =
\max_{j,k} \{ -d_{1k}, -p_j - e_j \}. \leqno(**)
$$
According to Lemma \ref{(c-1)-presentation-of-ideal} $P^*$ is a minimal
$(c-1)$-syzygy of $K_A(n+1) \cong \ER^{c-2}(P,R)$. Thus we obtain
$$
a(E^*) = a(P^*) \geq a(K_A(n+1)) + c - 1 = - e(X) - n - 1 + c - 1.
$$
Hence $(**)$ provides
$$
\max_{j,k} \{ d_{1k}, p_j + e_j \} \leq e(X) + n + 2 -c.
$$
This proves the estimate on the right-hand side of claim (a) and also of
claim (b) taking $(*)$ into account again.

Since $P$ does not split a free direct summand by Lemma
\ref{(c-1)-presentation-of-ideal} we obtain using Theorem I.4.1 of
\cite{hab} that
$$
\min_{j,k} \{ d_{2k} \} = a(P) \geq 1 - e^+(P^*) = 1 - e^+(E*).
$$
Hence $(**)$ implies the estimate on the left-hand side in claim (b).

It remains to show the estimate on the left-hand side in claim (a). Let
$C$ be the curve arising as intersection of $X$ with a general linear
subspace of dimension $c+1$. Then we get using for example Lemma
\ref{Omeag-res-unter-hyperpl-sect} that
$$
H^1_*(\cJ_C) \cong \bigoplus_{i = 0}^{n-1-c}
(H^{1+i}_*(\cJ_X)(-i))^{\binom{n-c-1}{i}}.
$$
It follows
$$
a(\HH^1(\cJ_C)) = \min \{ a(H^{1+i}_*(\cJ_X)) + i \} = \min \{ e_j + p_j \}
-1.
$$
Let $H$ be a general hyperplane and consider the exact sequence
$$
H_*^0(\cJ_{C \cap H}) \to H_*^1(\cJ_C)(-1) \stackrel{H}{\longrightarrow}
H_*^1(\cJ_C).
$$
Since $C$ is arithmetically Buchsbaum the map on the right-hand side is
zero. Thus we obtain
$$
a(H^1_*(\cJ_C)) \geq 0
$$
completing the proof in case $c \geq 3$.

Let now $c=2$. Then the claims follow similarly using the exact sequence
$$
0 \to R \to E^* \to P^* \to \ER^1(I,R) \to 0.
$$
\end{proof}

\begin{remark} In case $\codim X = 2$ there is a  considerably stronger
  result. In fact, Chang \cite{Chang-Diff-Geom}, Theorem 2.3 could
  characterize the possible twists in an $\Omega$-resolution precisely if
  $c = 2$. (Her results are also true in positive characteristic due to
  Walter \cite{Walter-transversality}.)  However, it seems to be rather
  difficult to achieve a similarly
  precise characterization  in higher
  codimension.  The corresponding problem is even open for arithmetically
  Cohen-Macaulay subschemes.
\end{remark}

\begin{example} Let $n \geq 3$ and $r \geq n-2$ be  two integers. There is
  an arithmetically Buchsbaum subscheme $X$ with $\Omega$-resolution
$$
0 \to (\cO_{\PP^n}(-r-1))^{n-2} \oplus \cO_{\PP^n}(n-2r-3) \to \cOP(-r+1)
\to \cJ_X \to 0
$$
where $e(X) = 2 (r-n+1)$ (cf.\ \cite{CGN2}, Lemma 6.10). For $X$ both
estimates in Proposition \ref{shifts-in-Omega-res}(b) are optimal whereas
the estimate on the right-hand side in Proposition
\ref{shifts-in-Omega-res}(a) is attained if $r = n-2$ and the one on the
left-hand side is optimal too if $r = n-2 = 1$.
\end{example}

\begin{remark} \label{laison-phenomena} It has been shown as Corollary 5.4
  in \cite{N-gorliaison} that an arithmetically Buchsbaum
  scheme $X \subset \PP^n$ is minimal in its even liaison class if the
  estimate on the
  right-hand side in Proposition
\ref{shifts-in-Omega-res}(a) is attained. Fixing any codimension $c \geq 2$
various such schemes do exist (cf.\  \cite{N-gorliaison}, Example 8.4).
\end{remark}

We conclude this section by pointing out that the
previous result implies \cite{Hoa-M}, Corollary 2.8. Recall that the
Castelnuovo-Mumford regularity of a subscheme $X$ (cf.\ \cite{M}) is the
integer
$$
\reg X = \max \{ i+1+ e(H^i_*(\cJ_X) \s i \geq 1 \}.
$$

\begin{corollary} If $X$ is arithmetically Buchsbaum of codimension $c$
  then
$$
e(X) + n + 2 - c \leq \reg X \leq e(X) + n + 3 - c.
$$
\end{corollary}

\begin{proof} By definition of the regularity we have
$$
\reg X - 1 = \max \{ i + e(H^i_*(\cJ_X)) \s i > 0 \} \geq e(X) + n + 1 - c.
$$
Using the notation of the previous proposition we also observe that
$$
\max \{ i + e(H^i_*(\cJ_X)) \s 1 \leq i \leq n-c \} = \max \{ e_j + p_j
\}.
$$
Thus the assertion follows by Proposition \ref{shifts-in-Omega-res} (a).
\end{proof}

In general it is not easy to compute the index of speciality. In
\cite{NS2}, Lemma 4.6 the estimate
$$
e(X) \leq  \deg X - \dim X - 2
$$
is shown. If $X$ is an integral subscheme there is the better estimate
(cf. \cite{NS2}, Lemma 4.6)
$$
e(X) < \left \lceil \frac{\deg X - 1}{c} \right \rceil - \dim X.
$$

\begin{remark}
(i)  We have left open the rather difficult smoothness and
  integrality questions. However, there are useful criteria for
  arithmetically Buchsbaum subschemes of codimension two due to Chang
  \cite{Chang-Diff-Geom}, Theorem 2.2 and 2.3. The smooth integral
  arithmetically Buchsbaum subschemes of arbitrary codimension which are
  divisors on a variety of
  minimal degree have been described in
  \cite{N-Buchsbaum-divisors}.

(ii) A different approach to study an arithmetically Buchsbaum
subscheme $X$
of $\PP^n$ has been proposed by M.\ Amasaki. Using the degrees of
the elements of a Gr\"obner basis of $I( X)$  in generic
coordinates he obtains a rough classification of arithmetically
Buchsbaum subschemes of projective space (cf.\ \cite{Amasaki-TAMS}
Theorem 4.5 and Corollary 6.7).

(iii) Notice that our structural approach allows to consider not
only subschemes of $\PP^n$ but even subschemes of some
arithmetically Gorenstein subscheme $G$. Indeed, our
characterization of
arithmetically Buchsbaum subschemes of $\PP^n$ has been obtained
by specializing the results about surjective-Buchsbaum subschemes
of $G$.

\end{remark}

\section{Quasi-Buchsbaum subschemes}
\label{quasi-Buchsbaum_schemes}

In this section we  show that every ideal of a Gorenstein ring
admits an exact sequence where all the occurring modules except the ideal
itself are Eilenberg-MacLane modules. As a consequence we obtain that every
equidimensional Cohen-Macaulay subscheme of projective space admits a locally
free resolution consisting of Eilenberg-MacLane bundles.
From this result we derive our characterization of  quasi-Buchsbaum
subschemes by means of weak $\Omega$-resolutions.

The construction in the proof of the following statement works for arbitrary
modules. However, for simplicity and keeping our applications in mind  we
restrict ourselves to saturated ideals.

\begin{theorem} \label{resolution_by_E-ML-modules} Let $I \subset R$ be
  a saturated ideal of  codimension $c$. Put $A =
  R/I$ and
$$s =  \min \{k \geq c \s \HH^j(A) = 0 \; \mbox{for all} \; j \; \mbox{with} \;
k+1 \leq j \leq n-k \}.
$$
Then $I$ admits an exact sequence
$$
0 \to E_s \to E_{s-1} \to \ldots \to E_1 \to I \to 0 \leqno(+)
$$
where all the modules $E_i$ are Eilenberg-MacLane modules of depth $2 i$ or
maximal Cohen-Macaulay modules, the module $E_i$  has finite  finite
projective dimension if $2 \leq i \leq s$ and $E_1$  has finite
finite projective dimension if and only if $I$ does.

Furthermore, the $R$-module  $E_i$
is free if $i >  n-c$ and we have
$$
\HH^{2i}(E_i) \cong \HH^i(A) \quad \mif i \leq \min \{n-c, s \}.
$$
Moreover, the sequence $(+)$ is uniquely determined up to isomorphisms of
exact sequences if it is impossible to cancel out free direct summands. In
this case $(+)$ is called {\rm minimal}. \\
If $s = n-c$ then all the modules $E_i$ have
cohomology of finite length if and only if $A$ is equidimensional and
locally Cohen-Macaulay.
\end{theorem}

\begin{proof} First we observe that the integer $s$ is
  well-defined. Indeed, the condition in its definition is trivially
  satisfied if $k \geq \frac{n}{2}$. Thus
we have $s = c$ if $c \geq \frac{n}{2}$ and
$c \leq s \leq \frac{n}{2}$ otherwise.  Let
$d = n+1-c$ denote the
  dimension of $R/I$. Now  we proceed in several steps.\\
\smallskip

\noindent
(I) We show that for any integer $k$ where $2 \leq k \leq \min \{d,
  c \}$ there is an exact sequence
$$
0 \to P_k \to E_{k-1} \to \ldots \to E_1 \to I \to 0 \leqno(C_k)
$$
where the modules $E_i$ have the properties as claimed  in our assertion
and $P_k$ is a module with  finite projective dimension and
$$
\HH^j(P_k) \cong \left \{ \begin{array}{ll}
0 & \mif j < 2k \\
\HH^{j-k}(A) & \mif 2k \leq j \leq n.
\end{array} \right.
$$
In order to prove this we induct on $k \geq 2.$ If $k = 2$ we take as
$(C_2)$ a minimal $(n-2)$-presentation of $I$. It has the desired properties
due to Lemma \ref{distribute_cohomology}. Now let $2 < k \leq \min \{d, c \}$. By induction we have an exact sequence
$$
0 \to P_{k-1} \to E_{k-2} \to \ldots \to E_1 \to I \to 0 \leqno(C_{k-1})
$$
where $\HH^{d+k-1}(P_{k-1}) \cong \HH^d(A) \neq 0$ because $2k-2 < d + k-1 <
d+c = n+1$.
Consider a minimal $(n+2-2k)$-presentation
$$
0 \to P_k \to E_{k-1} \to P_{k-1} \to 0
$$
of the module $P_{k-1}$. Then we obtain by induction and Lemma
\ref{distribute_cohomology}
$$
\HH^{j}(E_{k-1}) \cong \left \{
  \begin{array}{@{\hspace{.0cm}}l@{\hspace{.2cm}}l} 
\HH^j(P_{k-1}) & \mif j < 2k-1 \\
0 & \mif 2k-1 \leq j \leq n
\end{array} \right.
\cong \left \{ \begin{array}{@{\hspace{.0cm}}l@{\hspace{.2cm}}l}
\HH^{k-1}(A) & \mif j = 2k-2 \\
0 & \mif j \neq 2k-2,n+1
\end{array} \right.
$$
and
$$
\HH^j(P_k) \cong \left \{ \begin{array}{ll}
0 & \mif j \leq 2k-1 \\
\HH^{j-1}(P_{k-1}) & \mif 2k \leq j \leq n
\end{array} \right.
\cong \left \{ \begin{array}{ll}
0 & \mif j < 2k \\
\HH^{j-k}(A) & \mif 2k \leq j \leq n.
\end{array} \right.
$$
In particular, because of $2k-2  < d + k-1 \leq n$  we see that
$\HH^{d+k-1}(E_{k-1}) =0$ is not isomorphic
to  $\HH^{d+k-1}(P_{k-1}) \neq 0$. Thus, $P_k$ cannot be trivial.
Therefore splicing together the sequence $(C_{k-1})$ and the minimal
$(n+2-2k)$-presentation of $P_{k-1}$ yields the desired exact sequence
$(C_{k})$. \\
\smallskip

\noindent
(II) Now we distinguish two cases. \\
{\it Case 1}: Suppose that $d < c$. Thus we get $s = c$. Let $d \geq
2$. Then we have an exact
sequence
$(C_{d})$ by step (I) of the proof. The module $P_{d}$ has finite
projective dimension and satisfies
$$
\HH^j(P_{d}) \cong \left \{ \begin{array}{ll}
0 & \mif j \leq 2d-1 \\
\HH^{j-d}(A) & \mif  2d \leq j \leq n.
\end{array} \right.
$$
Thus $P_{d}$ has depth $2d$ since $2d \leq n$ due to our assumption
of this case. It follows by the
Auslander-Buchsbaum formula that $P_{d}$ has a minimal free resolution
$$
0 \to F_{n+1-2d} \to \ldots \to F_0 \to P_d \to 0.
$$
Splicing together this resolution and the sequence $(C_d)$  yields the
desired sequence $(+)$ of length $s = c$ where we relabel the modules $F_i$
as $E_{d+i}$.

If $d = 1$ then we consider a minimal $(n-1)$-presentation of $I$
$$
0 \to P \to E_1 \to I \to 0.
$$
Since $I$ is saturated we obtain by Lemma \ref{distribute_cohomology} that
$E_1$ is a maximal Cohen-Macaulay module and that $P$ is an
Eilenberg-MacLane module of depth $3$. Thus, replacing $P$ by its minimal
free resolution provides the desired sequence $(+)$ of length $s = n$. \\
{\it Case 2}: Suppose that $d \geq c$.
 Thus we have an exact sequence $(C_{c})$ by step
(I). We claim that $I$  admits even an exact sequence $(C_s)$ with the
properties described in step (I). Indeed, continuing in the fashion of
step (I) we just need to check that the modules $P_c,\ldots,P_s$ are
non-trivial. If $t =c $ we are done by step (I). If $c <k \leq  s$
and we
have constructed $P_{k-1}$ we know by the definition of $s$ that there is an
integer $j$ such that $2k-1 \leq  j \leq n$ and $\HH^j(P_{k-1}) \cong
\HH^{j-k+1}(A) \neq 0$.  Since $\HH^j(E_{k-1}) = 0$ if $2k-1 \leq  j \leq
n$ it follows that the module $P_k$ cannot be trivial by considering the
minimal  $(n+2-2k)$-presentation
$$
0 \to P_k \to E_{k-1} \to P_{k-1} \to 0.
$$

Now we turn our attention to the module $P_s$ occurring in  $(C_{s})$.  It
has depth $\geq 2s$ by construction. Furthermore, we have  that $\HH^j(P_s) \cong
\HH^{j-s}(A) \neq 0$ if $2s+1 \leq j \leq n$ due to  the definition of
$s$. Therefore $P_s$
is an Eilenberg-MacLane module. Thus,  defining $E_s = P_s$ the sequence
$(C_{s})$ gives the desired sequence $(+)$. \\
\smallskip

\noindent
(III) Given a minimal sequence $(+)$ we can reverse the constructions in
steps (I) and (II) by shopping it into  shorter exact sequences. Thus the
claimed uniqueness is a consequence of the uniqueness properties of
$q$-presentations and free resolutions.

The final claim follows because  in case $s \geq n-c$ all the modules
$E_i$  in $(+)$  have
cohomology of finite length if and only if  $A$ has cohomology of finite
length. The latter
is true if and only if $A$ is equidimensional and locally Cohen-Macaulay.
\end{proof}

\begin{remark} Using the notation of the proposition above let us assume
  that $I$ is a perfect ideal. Then we obtain $s = c$ and the minimal
  sequence $(+)$ is nothing else than the minimal free resolution of $I$.
\end{remark}

Observe that the sequence $(+)$ above can have a length which is less than
the dimension of
$A$. This is also reflected in the following result where we use the
modules $G_i$ defined immediately before Lemma~\ref{properties-of-indecomposable-max-Buchs}.

\begin{corollary} \label{half-quasi-BM}
Let $I \subset R$ be a saturated ideal of finite projective dimension. Put
$s = \min \{k \geq c \s \HH^j(R/I) = 0 \; \mbox{for all} \; j \; \mbox{with} \;
k+1 \leq j \leq n-k \}$ and $v = \min \{n-c, s\}$. Then the following
conditions are equivalent:
\begin{itemize}
\item[(a)] It holds
$$
\fm \cdot \HH^i(R/I)= 0 \quad \mif  1 \leq i \leq v.
$$
\item[(b)]  $I$ admits an exact sequence
\begin{multline*}
0 \to F_s \to \ldots \to F_{v+1} \to F_v \oplus \bigoplus_j
(G_{2v}(-e_{v.j}))^{s_{v.j}} \to \ldots \\
 \to F_1
\oplus \bigoplus_j (G_{2}(-e_{1.j}))^{s_{1.j}} \to I \to 0
\end{multline*}
where the integers $s_{i.j}$ are non-negative.
\end{itemize}
\end{corollary}

\begin{proof} Consider the exact sequence $(+)$ given by Theorem
  \ref{resolution_by_E-ML-modules}. If (a) is satisfied then the modules
  $E_i$ are free if $v < i \leq s$ and are Eilenberg-MacLane modules of finite projective dimension where
$\fm \cdot \HH^{2i}(E_i) = 0$ otherwise. Hence the $E_i, 1 \leq i \leq v,$
are surjective-Buchsbaum modules by Lemma
\ref{hinreichende-Bedingung-fuer-surj-Buchs}. Therefore Proposition
\ref{char-of-surj-Buchsbaum} shows $E_i \cong F_i \oplus \bigoplus_j
(G_{2i}(-e_{i.j}))^{s_{i.j}}$ for some free module $F_i$. Thus the sequence
$(+)$ is of the form as claimed in (b).

In order to show the converse we use the sequence in (b) to compute
$$
\HH^i(R/I) \cong \HH^{2i}(G_{2i}) \quad \mif 1 \leq i \leq v.
$$
Since the modules $G_i$ are Buchsbaum we obtain (a).
\end{proof}

Now we specialize our results to the case where $R$ is a polynomial
ring.

\begin{corollary} \label{resolution_by_Eil-MacLane-bundles} Let $X \subset \PP^n$ be an equidimensional
  Cohen-Macaulay subscheme of codimension $c$. Then $X$ admits a locally
  free resolution
$$
0 \to \cE_s \to \ldots \to \cE_1 \to \cJ_X \to 0 \leqno(*)
$$
where $s = \min \{k \geq c \s H^j_*(\cJ_X) = 0 \; \mbox{for all} \; j \;
\mbox{with} \; k+1 \leq j \leq n-k \}$ and $\cE_1,\ldots,\cE_s$ are
Eilenberg-MacLane bundles which split as a direct sum of line bundles
unless $\cE_i$ has depth $2i - 1$ and  $1 \leq i \leq v = \min \{
n-c, s\}$.

For any such resolution it holds
$$
H^i_*(\cJ_X) \cong H^{2i-1}_*(\cE_i) \quad \mif 1 \leq i \leq v.
$$
\end{corollary}

\begin{proof} Consider the sequence $(+)$ of the ideal $I_X = H^0_*(\cJ_X)$
  provided by Theorem  \ref{resolution_by_E-ML-modules}. Due to the assumption on $X$ the modules $E_1,\ldots,E_s$
  have cohomology of finite length. Thus the sheaves $\cE_i = \tilde{E_i}$
  are locally free and we obtain the sequence $(*)$ as sheafification of
  $(+)$.

Conversely, taking global sections in any sequence $(*)$ we get a sequence
of the form $(+)$ as in Theorem  \ref{resolution_by_E-ML-modules}. Thus
the uniqueness properties of $(+)$ show that the sequence $(*)$ is indeed a
locally free resolution.

The final assertion follows by tracing the cohomology along short exact
sequences.
\end{proof}

We want to give a name to a certain type of exact  sequences.

\begin{definition} \label{def-weak-Omega-resolution} A
  subscheme $X \subset \PP^n$  of codimension $c$ is said to have a {\it
    weak $\Omega$-resolution} if there exists an exact sequence
\begin{multline*}
0  \longrightarrow \cF_s \stackrel{\alpha_s}{\longrightarrow}  \ldots
\stackrel{\alpha_{v+2}}{\longrightarrow} \cF_{v+1} \longrightarrow
\cF_v \oplus  \bigoplus_j
(\cOP^{2v-1}(-e_{v.j}))^{s_{v.j}} \stackrel{\alpha_v}{\longrightarrow} \\
\ldots \stackrel{\alpha_{2}}{\longrightarrow} \cF_1
\oplus \bigoplus_j (\cOP^{1}(-e_{1.j}))^{s_{1.j}} \longrightarrow I
\longrightarrow 0
\end{multline*}
where $s = \min \{k \geq c \s H^j_*(\cJ_X)
  = 0 \; \mbox{for all} \; j \; \mbox{with} \; k+1 \leq j \leq n-k \}$, $v = \min\{n-c, s\}$, $\cF_1,\ldots,\cF_s$ are direct sums
of line bundles and for fixed $i$ the integers $e_{i.j}$ are pairwise
distinct.

The weak $\Omega$-resolution is called {\it minimal} if there is no line
bundle $\cL$ in the sequence such that the restriction of $\alpha_i \  (2 \leq
i \leq s)$ to $\cL$ induces an isomorphism of $\cL$ onto $\cL$.
\end{definition}

Subschemes having a weak $\Omega$-resolution can be characterized
cohomologically.

\begin{theorem} \label{resolution-by-Omegas} Let $X \subset \PP^n$ be a
  subscheme of codimension $c$. Let $s = \min \{k \geq c \s H^j_*(\cJ_X)
  = 0 \; \mbox{for all} \; j \; \mbox{with} \; k+1 \leq j \leq n-k \}$ and
  $v = \min\{n-c, s\}$.
Then the following conditions are equivalent:
\begin{itemize}
\item[(a)] It holds
$$
\fm \cdot H^i_*(\cJ_X)= 0 \quad \mif  1 \leq i \leq v.
$$
\item[(b)]  $X$ admits a weak $\Omega$-resolution
\begin{multline*}
0 \to \cF_s \to \ldots \to \cF_{v+1} \to \cF_v \oplus \bigoplus_j
(\cOP^{2v-1}(-e_{v.j}))^{s_{v.j}} \to \\
\ldots \to \cF_1
\oplus \bigoplus_j (\cOP^{1}(-e_{1.j}))^{s_{1.j}} \to I \to 0
\end{multline*}
where  $s_{i.j} = h^i(\cJ_X(j))$.
\end{itemize}
Moreover, the weak $\Omega$-resolution is a locally free resolution, i.e.,
the minimal weak $\Omega$-resolution of $X$ is uniquely determined.
\end{theorem}

\begin{proof} Recall that $\widetilde{G_{i+1}} \cong \cOP^i$ if $R = K[x_0,\ldots,x_n]$ is a
  polynomial ring. Thus, the result follows by the Corollaries
  \ref{half-quasi-BM} and
  \ref{resolution_by_Eil-MacLane-bundles}
\end{proof}

\begin{remark}
In Theorem \ref{intro_char_quasi-BM} of the introduction we gave a
characterization of
quasi-Buchsbaum subschemes  $X \subset \PP^n$ where $c \geq \frac{n}{2}$.
This result is
is covered by Theorem \ref{resolution-by-Omegas} because $s = c$ and $v =
n-c$ if $c \geq \frac{n}{2}$.
\end{remark}

\section{Surfaces in $\PP^4$} \label{section_surfaces}

Now we want to specialize some results to surfaces of codimension two. This
gives a new prospect on the construction of smooth surfaces in
$\PP^4$. Then we derive a criterion on a weak $\Omega$-resolution of a
surface $X$ in order to obtain a characterization of the arithmetically
Buchsbaum surfaces among the quasi-Buchsbaum surfaces. This criterion could
be extended to arbitrary quasi-Buchsbaum subschemes $X \subset \PP^n$ such
that $\dim X \leq \codim X$. We won't pursue this here. Instead, we explain
how the criterion can be applied.

\begin{proposition} \label{Eil-MacLane_res_of_surfaces} Let $X\subset
  \PP^4$ be a $2$-dimensional subscheme. Then $X$ has a  presentation
$$
0 \to \cE_2 \stackrel{\ffi}{\longrightarrow} \cE_1 \to \cJ_X \to 0
$$
where $\cE_2$ and $\cE_1$ are Eilenberg-MacLane sheaves such that
$H^i_*(\cE_1) = 0$ if $i = 2, 3$ and $H^i_*(\cE_2) = 0$ if $i = 1,
2$.

If $\ffi$ is a minimal morphism, i.e.\ it does not map a line bundle
summand of $\cE_2$ onto a line bundle summand of $\cE_1$, then the
presentation is uniquely determined (up to isomorphisms of exact
sequences).

Furthermore, the sheaf $\cE_1$ is torsion-free, $\cE_2$ is
reflexive and
$$
H^1_*(\cE_1) \cong H^1_*(\cJ_X) \; \mbox{and} \; H^2_*(\cE_2) \cong
H^2_*(\cJ_X).
$$

Moreover, we have:
\begin{itemize}
\item[(a)] $X$ is equidimensional if and only if $\cE_1$ is a vector bundle
  and $\dim H^3_*(\cE_2)^{\vee} \leq 1$.
\item[(b)] $X$ is equidimensional and Cohen-Macaulay if and only if $\cE_1$
  and $\cE_2$  are vector bundles.
\end{itemize}
\end{proposition}

\begin{proof} The existence of the asserted presentation is a consequence
  of Theorem \ref{resolution_by_E-ML-modules}. Since $\dim
  H^1_*(\cJ_X)^{\vee} \leq 1$ and $\dim H^2_*(\cJ_X)^{\vee} \leq 2$,
  Proposition \ref{k-syzch} shows that $\cE_1$ is torsion-free and that
  $\cE_2$ is reflexive.

Claim (b) is a special case of Corollary
\ref{resolution_by_Eil-MacLane-bundles}. Claim (a) follows by combining
Proposition \ref{k-syzch}
and the cohomological characterization of
  equidimensionality \cite{N-gorliaison}, Lemma 2.11.
\end{proof}

\begin{remark} \label{remark_surface_construction}
(i) Using Horrocks' \cite{Ho1} characterization of stable equivalence
classes of vector bundles on $\PP^4$ Bolondi \cite{Bolondi-surfaces} also
obtained the presentation of the statement in the case of equidimensional
Cohen-Macaulay surfaces. But he did not establish its uniqueness
property.

(ii) Using the Eagon-Northcott complex we may think of our subscheme $X$ as
the degeneracy locus of the map $\ffi$. Note that the idea of
constructing surfaces as degeneracy locus of vector bundles has been
systematically exploited in \cite{DES-surfaces}. There the bundles are
chosen in accordance with Beilinson's spectral sequence and not
necessarily Eilenberg-MacLane sheaves.

(iii) The presentation above allows to construct an equidimensional
surface with prescribed cohomology. Indeed, the cohomology of $X$ determines the
cohomology of the Eilenberg-MacLane sheaves $\cE_1$ and $\cE_2$ which in
turn determine $\cE_1$ and $\cE_2$ itself up to direct sums of line bundles
according to Proposition \ref{elc}. Then  we have to choose a sufficiently
general map $\ffi: \cE_2 \to \cE_1$. If the degeneracy locus has
codimension two then it is the desired surface.
\end{remark}

Now we want to compare quasi-Buchsbaum and arithmetically Buchsbaum
surfaces.
Let $X \subset \PP^4$ be a quasi-Buchsbaum surface. According to Theorem
\ref{resolution-by-Omegas} it has a weak $\Omega$-resolution
$$
0 \to \cF_2 \oplus \bigoplus_j (\Omega_{\PP^4}^{3}(-e_{2.j}))^{s_{2.j}}
\stackrel{\ffi}{\longrightarrow}
 \cF_1 \oplus \bigoplus_j (\Omega_{\PP^4}^{1}(-e_{1.j}))^{s_{1.j}} \to \cJ_X
 \to 0
$$
where $\cF_1, \cF_2$ are direct sums of line bundles on $\PP^4$.  Let $\cP$
be a direct sum of line bundles such that there is an epimorphism $\delta:
\cP \to \cF_1 \oplus \bigoplus_j
(\cO_{\PP^4}^{1}(-e_{1.j}))^{s_{1.j}}$. Using this notation we can
distinguish between quasi-Buchsbaum and arithmetically Buchsbaum
subschemes  as follows.

\begin{proposition} \label{aBM-among-quasi-BM} The surface $X$ is
  arithmetically Buchsbaum if and only if there is a morphism
  $\alpha: \cF_2
  \oplus \bigoplus_j (\Omega_{\PP^4}^{3}(-e_{2.j}))^{s_{2.j}} \to \cP$ such
  that the following diagram is commutative
$$
\xy\xymatrixrowsep{1.5pc}\xymatrixcolsep{1.4pc}\xymatrix{
0 \ar @{->}[r] &  \cF_2  \oplus  
\bigoplus_j (\Omega_{\PP^4}^{3}(-e_{2.j}))^{s_{2.j}}   
\ar @{->}[dr]_-{\alpha}   \ar
@{->}^-{\ffi}[r] & \cF_1 \oplus \bigoplus_j
(\Omega_{\PP^4}^{1}(-e_{1.j}))^{s_{1.j}}  \ar
 @{->}[r] & \cJ_X \ar @{->}[r] &  0 \\
&   &  \cP. \ar @{->}[u]_-{\delta}
}
\endxy
$$
\end{proposition}

\begin{proof} Assume that such a map $\alpha$ exists. Consider the
  exact commutative diagram
$$
\begin{array}{c@{\hspace{0.2cm}}c@{\hspace{0.2cm}}c@{\hspace{0.2cm}}c@{\hspace{0.2cm}}c@{\hspace{0.2cm}}c@{\hspace{0.2cm}}c@{\hspace{0.2cm}}c@{\hspace{0.2cm}}c}   
& & & 0 \\
& & & \Big\uparrow \\
0 \to & \cF_2 \oplus \bigoplus_j (\Omega_{\PP^4}^{3}(-e_{2.j}))^{s_{2.j}} &
\stackrel{\ffi}{\longrightarrow} &
 \cF_1 \oplus \bigoplus_j (\Omega_{\PP^4}^{1}(-e_{1.j}))^{s_{1.j}} &
 \longrightarrow & 
 \cJ_X & \to 0 \\
&   &  & \Big\uparrow\!{\scriptstyle \delta}  & & \Big\uparrow\!{\scriptstyle =} & \\
0 \to & \ker \gamma & \stackrel{\varepsilon}{\longrightarrow}  & \cP &
\stackrel{\gamma}{\longrightarrow} & \cJ_X & \to 0 \\
& & & \Big\uparrow & \\[1ex]
& & & \cP' \oplus \bigoplus_j (\Omega_{\PP^4}^{2}(-e_{1.j}))^{s_{1.j}} & \\[1ex]
& & & \Big\uparrow \\
& & & 0 \\
\end{array}
$$
where $\cP'$ is a sum of line bundles. The Snake lemma provides the exact
sequence
$$
0 \to \cP' \oplus \bigoplus_j (\Omega_{\PP^4}^{2}(-e_{1.j}))^{s_{1.j}} \to
\ker \gamma \stackrel{\beta}{\longrightarrow} \cF_2 \oplus \bigoplus_j
(\Omega_{\PP^4}^{3}(-e_{2.j}))^{s_{2.j}} \to 0.
$$
This sequence splits because $\varepsilon^{-1} \circ \alpha$ is a splitting
map. Indeed, our assumption implies $\im \alpha \subset \ker \gamma = \im
\vep$ and $\beta \circ (\varepsilon^{-1} \circ \alpha)$ is the identity map
on $\cF_2 \oplus \bigoplus_j (\Omega_{\PP^4}^{3}(-e_{2.j}))^{s_{2.j}}$.
Therefore the second row of the first diagram implies that $X$ is
arithmetically Buchsbaum according to Corollary
\ref{geom-char-of-aBM-schemes}.

Now let us assume that $X$ is arithmetically Buchsbaum. Then $X$ has a
locally free resolution
$$ 0 \to \cF_2 \stackrel{\psi}{\longrightarrow} \cF_1' \oplus
\bigoplus_j (\Omega_{\PP^4}^{2}(-e_{2.j}))^{s_{2.j}}
\oplus \bigoplus_j (\Omega_{\PP^4}^{1}(-e_{1.j}))^{s_{1.j}} \to \cJ_X \to 0
$$
according to Corollary
\ref{geom-char-of-aBM-schemes}.

Resolving $\bigoplus_j (\Omega_{\PP^4}^{2}(-e_{2.j}))^{s_{2 .j}}$
and using the Snake lemma we get an exact commutative diagram
$$
\begin{array}{c@{\hspace{0.2cm}}c@{\hspace{0.2cm}}c@{\hspace{0.2cm}}c@{\hspace{0.2cm}}c@{\hspace{0.2cm}}c@{\hspace{0.2cm}}c@{\hspace{0.2cm}}c@{\hspace{0.2cm}}c} 
& 0 & & 0 \\
&  \Big\uparrow & & \Big\uparrow \\[1ex]
0 \to & \cF_2  &
\stackrel{\psi}{\longrightarrow} & \begin{array}{c}
 \cF_1' \oplus \bigoplus_j
 (\Omega_{\PP^4}^{1}(-e_{1.j}))^{s_{1.j}} \\
\oplus \\
\bigoplus_j (\Omega_{\PP^4}^{2}(-e_{2.j}))^{s_{2.j}}
\end{array}
  & {\longrightarrow}  &
 \cJ_X & \to 0 \\[1ex]
& \Big\uparrow  &  & \Big\uparrow\!{\scriptstyle }  & &
\Big\uparrow\!{\scriptstyle =} & \\
0 \to &
\ker \gamma &
{\longrightarrow}  & \cF_1 \oplus \bigoplus_j
 (\Omega_{\PP^4}^{1}(-e_{1.j}))^{s_{1.j}} &
\stackrel{\gamma}{\longrightarrow} & \cJ_X & \to 0 \\
& \Big\uparrow & & \Big\uparrow & \\
& \bigoplus_j (\Omega_{\PP^4}^{3}(-e_{2.j}))^{s_{2.j}} & =  &  \bigoplus_j
  (\Omega_{\PP^4}^{3}(-e_{2.j}))^{s_{2.j}} & \\
&  \Big\uparrow & & \Big\uparrow \\
& 0 & & 0 \\
\end{array}
$$
where $\cF_1$ is a direct sum of line bundles.  Since $\cF_2$ is a ´direct
sum of line bundles the vertical sequence on the left-hand side
must split. Repeating the argument where we also resolve $\bigoplus_j
(\Omega_{\PP^4}^{1}(-e_{1.j}))^{s_{1.j}}$  we
get an exact commutative diagram
$$
\begin{array}{c@{\hspace{0.2cm}}c@{\hspace{0.2cm}}c@{\hspace{0.2cm}}c@{\hspace{0.2cm}}c@{\hspace{0.2cm}}c@{\hspace{0.2cm}}c@{\hspace{0.2cm}}c@{\hspace{0.2cm}}c}
0 \to & \cF_2 \oplus \bigoplus_j (\Omega_{\PP^4}^{3}(-e_{2.j}))^{s_{2.j}} &
\stackrel{\psi'}{\longrightarrow} &
 \cF_1 \oplus \bigoplus_j (\Omega_{\PP^4}^{1}(-e_{1.j}))^{s_{1.j}} & \to
 \cJ_X  \to 0 \\[1ex]
 & \Big\uparrow{\scriptstyle \pi}  &  & \Big\uparrow\!{\scriptstyle \delta'}  &
\Big\uparrow\!{\scriptstyle =} & \\ [1ex]
0 \to &
\begin{array}{c}
\cF_2 \oplus \bigoplus_j (\Omega_{\PP^4}^{3}(-e_{2.j}))^{s_{2.j}}
\\
\oplus \\
\bigoplus_j (\Omega_{\PP^4}^{2}(-e_{1.j}))^{s_{1.j}}
\end{array}
 &
\stackrel{\psi''}{\longrightarrow} &
 \cF_1 \oplus \cP' & \to
 \cJ_X  \to 0
\end{array}
$$
where $\pi$ is the canonical projection and $\cP'$ a direct sum of line bundles.
Hence, denoting by $\alpha$ the
restriction of $\psi''$ to  $\cF_2 \oplus \bigoplus_j
(\Omega_{\PP^4}^{3}(-e_{2.j}))^{s_{2.j}}$
we obtain the commutative diagram
$$
\xy\xymatrixrowsep{1.5pc}\xymatrixcolsep{1.4pc}\xymatrix{
0 \ar @{->}[r] & \cF_2 \oplus \bigoplus_j
(\Omega_{\PP^4}^{3}(-e_{2.j}))^{s_{2.j}} \ar @{->}[dr]_-{\alpha}   \ar
@{->}^-{\psi'}[r] & \cF_1 \oplus \bigoplus_j
(\Omega_{\PP^4}^{1}(-e_{1.j}))^{s_{1.j}}  \ar
 @{->}[r] & \cJ_X \ar @{->}[r] &  0 \\
&   & \cF_1 \oplus \cP'. \ar @{->}[u]_-{\delta'}
}
\endxy
$$
This proves our claim for this particular weak
$\Omega$-resolution and the particular epimorphism $\delta'$. But splitting
off redundant line bundles we get  a minimal weak
$\Omega$-resolution and a minimal epimorphism $\delta$ and thus our claim
holds in this situation, too. Due to the
uniqueness of minimal $\Omega$-resolutions and  minimal free resolutions we
can conclude that {\it any} weak $\Omega$-resolution admits a morphism
$\alpha$ as claimed.
\end{proof}

Now we want to show how the criterion above can be applied.

\begin{example} \label{example-B1.15} Every surface $X \subset \PP^4$ with
  a weak $\Omega$-resolution
$$
0 \to (\Omega_{\PP^4}^{3}(-1))^2 \stackrel{\ffi}{\longrightarrow}
\cO_{\PP^4}(-4)  \oplus (\Omega_{\PP^4}^{1}(-3))^2 \to \cJ_X \to 0
$$
is quasi-Buchsbaum but not arithmetically Buchsbaum. Note that $X$ is a
smooth rational surface of degree $10$ if $\ffi$ is general enough (cf.\
\cite{DES-surfaces}, Example B1.15).

We need only to show that $X$ is  not arithmetically Buchsbaum. Assuming
the contrary there is  by Proposition \ref{aBM-among-quasi-BM} a 
commutative diagram
$$
\xy\xymatrixrowsep{1.5pc}\xymatrix{
0 \ar @{->}[r] &
(\Omega_{\PP^4}^{3}(-1))^{2} \ar @{->}[dr]_-{\alpha}   \ar
@{->}^-{\ffi}[r] & \cO_{\PP^4}(-4) \oplus
(\Omega_{\PP^4}^{1}(-3))^{2}  \ar
 @{->}[r] & \cJ_X \ar @{->}[r] &  0 \\
&   &  \cO_{\PP^4}(-4) \oplus (\cO_{\PP^4}(-5))^{20}. \ar @{->}[u]_-{\delta}
}
\endxy
$$
Since $Hom_{\cO_{\PP^4}}(\Omega_{\PP^4}^{3}(-1), \cO_{\PP^4}(-5)) \cong
H^0(\Omega_{\PP^4}^{1}(-1)) = 0$ we see that $\alpha$ must induce an
injective map $(\Omega_{\PP^4}^{3}(-1))^2 \to \cO_{\PP^4}(-4)$ which does
not exist. This contradiction proves our claim.
\end{example}

\begin{example} \label{example-B7.15} Every surface $X \subset \PP^4$ with
  a weak $\Omega$-resolution
$$
0 \to (\cO_{\PP^4}(-5))^2  \oplus \Omega_{\PP^4}^{3}(-1)
\stackrel{\ffi}{\longrightarrow}
(\cO_{\PP^4}(-4))^3  \oplus \Omega_{\PP^4}^{1}(-3) \to \cJ_X \to 0
$$
is quasi-Buchsbaum but not arithmetically Buchsbaum. Note that $X$ is a
 smooth elliptic surface of degree $10$ if $\ffi$ is general enough (cf.\
\cite{DES-surfaces}, Example B7.5).

The proof that $X$ is not arithmetically Buchsbaum is similar to the one
just given. Indeed, if $X$ were arithmetically Buchsbaum we had a map
$$ 
\alpha: \Omega_{\PP^4}^{3}(-1) \to (\cO_{\PP^4}(-4))^3  \oplus
(\cO_{\PP^4}(-5))^{10}
$$ 
 which must induce an injective map
$\Omega_{\PP^4}^{3}(-1) \to (\cO_{\PP^4}(-4))^3$. Again,  rank considerations
show that such an embedding cannot exist.
\end{example}


\end{document}